\documentclass[fleqn,10pt]{SelfArx} 
\usepackage[usenames,dvipsnames]{xcolor}
\usepackage[english]{babel} 
\usepackage[color=cyan]{todonotes}
\usepackage{lipsum} 
\usepackage{float}
\usepackage{breakcites}

\usepackage[titletoc,title]{appendix}
\usepackage[usenames,dvipsnames]{xcolor}
\usepackage[most]{tcolorbox}
\usepackage{tabularx}
\usepackage{array}
\usepackage{colortbl}
\tcbuselibrary{skins}
\usepackage{fixltx2e}
\newcolumntype{Y}{>{\raggedleft\arraybackslash}X}

\newcommand{\nocontentsline}[3]{}
\newcommand{\tocless}[2]{\bgroup\let\addcontentsline=\nocontentsline#1{#2}\egroup}

\newcommand\blfootnote[1]{%
  \begingroup
  \renewcommand\thefootnote{}\footnote{#1}%
  \addtocounter{footnote}{-1}%
  \endgroup
}

\usepackage{graphicx}
\usepackage{subcaption}

\tcbset{tab1/.style={fonttitle=\bfseries\large,fontupper=\normalsize\sffamily,
colback=yellow!10!white,colframe=red!75!black,colbacktitle=Salmon!40!white,
coltitle=black,center title,freelance,frame code={
\foreach \n in {north east,north west,south east,south west}
{\path [fill=red!75!black] (interior.\n) cycle (3mm); };},}}

\tcbset{tab2/.style={enhanced,fonttitle=\bfseries,fontupper=\normalsize\sffamily,
colback=yellow!10!white,colframe=red!50!black,colbacktitle=Salmon!40!white,
coltitle=black,center title}}

\usepackage{tcolorbox}
\usepackage{tabularx}
\usepackage{array}
\usepackage{colortbl}
\usepackage{breqn}
\tcbuselibrary{skins}

\usepackage{tablefootnote}

\newcolumntype{Y}{>{\raggedleft\arraybackslash}X}

\tcbset{tab1/.style={fonttitle=\bfseries\large,fontupper=\normalsize\sffamily,
colback=yellow!10!white,colframe=red!75!black,colbacktitle=Salmon!40!white,
coltitle=black,center title,freelance,frame code={
\foreach \n in {north east,north west,south east,south west}
{\path [fill=red!75!black] (interior.\n) cycle (3mm); };},}}

\tcbset{tab2/.style={enhanced,fonttitle=\bfseries,fontupper=\normalsize\sffamily,
colback=green!10!white,colframe=green!50!black,colbacktitle=green!40!white,
coltitle=black,center title}}

\definecolor{color1}{RGB}{0,139,0} 
\definecolor{color2}{RGB}{154,255,154} 

\numberwithin{equation}{section}

\newcommand{\Z}{\mathbb{Z}} 
\newcommand{\N}{\mathbb{N}}

\usepackage{hyperref} 
\usepackage{url}
\usepackage[framemethod=tikz]{mdframed}
\usepackage{cuted}
\usepackage{amsmath}
\usepackage{enumitem}
\hypersetup{hidelinks,colorlinks,breaklinks=true,urlcolor=color2,citecolor=color1,linkcolor=color1,bookmarksopen=false,pdftitle={Title},pdfauthor={Author}}
\usepackage{amssymb,amsmath,amsthm,amsfonts}
\newtheorem{theorem}{Theorem}
\newtheorem{definition}{Definition}
\JournalInfo{$~$} 
\Archive{} 

\PaperTitle{Vehicle Routing with Drones} 

\Authors{Rami Daknama\textsuperscript{1}*, Elisabeth Kraus\textsuperscript{1}*} 
\affiliation{\textsuperscript{1}\textit{Department of Mathematics, Ludwig-Maximilians-Universität München, Munich, Germany} \\ \{rami.daknama, elisabeth.kraus\}@math.lmu.de} 
\affiliation{* Both authors contributed equally to this work.}  

\newcommand{\keywordname}{Keywords} 

\Abstract{
We introduce a package service model where trucks as well as drones can deliver packages. Drones can travel on trucks or fly; but while flying, drones can only carry one package at a time and have to return to a truck to charge after each delivery.
We present a heuristic algorithm to solve the problem of finding a good schedule for all drones and trucks. The algorithm is based on two nested local searches, thus the definition of suitable neighbourhoods of solutions is crucial for the algorithm. Empirical tests show that our algorithm performs significantly better than a natural Greedy algorithm. Moreover, the savings compared to solutions without drones turn out to be substantial, suggesting that delivery systems might considerably benefit from using drones in addition to trucks.
}

\usepackage{tablefootnote}

\begin{document}

\flushbottom 

\maketitle 

\blfootnote{We used a \LaTeX \ template by Mathias Legrand (extensively modified by "Vel") which can be downloaded under http://www.latextemplates.com/template/stylish-article}

\tableofcontents 

\thispagestyle{empty} 

\section{Introduction}

Many enterprises which are interested in logistics have put considerable effort into developing delivery systems using drones for transporting packages, e.g. Amazon's project PrimeAir, that just made its first demo delivery in the U.S. in March 2017.
This gives rise to the following optimization problem. Imagine you have a fleet of trucks and you are given packages which have to be delivered to given positions in the plane. Moreover, there is a number of drones given, that can be carried by any truck on its roof, and each drone can carry one (arbitrary) package at a time while flying.

\blfootnote{We used a \LaTeX \ template by Mathias Legrand (extensively modified by "Vel") which can be downloaded under http://www.latextemplates.com/template/stylish-article}

Each time after having delivered a package, a drone returns to a truck and has to charge on the truck's roof by staying on it for at least one edge (meaning between two package delivery positions) of the truck tour.
Given this situation, a possible objective is to minimize the average time until a package is delivered.

  The problem generalizes the travelling salesman problem, thus it is NP-hard ---  indeed, in practice the problem seems much harder than the TSP.
 This complexity motivates the development of heuristics to solve the problem approximately. We develop a heuristic that uses two nested local searches. For implementing these heuristics, we exploit the metaheuristic framework JAMES \cite{james} which provides different local search methods when you customize initial solutions and neighbourhood definitions. Our contributions are the following:
\begin{itemize}
\item We introduce the problem (Vehicle Routing with Drones (VRD)) formally (different but related models already were introduced earlier, compare chapter \ref{literature}). 
\item We provide an equivalent characterization of feasibility of solutions to the VRD. 
\item We develop a heuristic algorithm to obtain solutions of good quality.
\item We implement the algorithm obtaining computational results; thereby we study the effects of different problem parameters.
\item We evaluate the computational results by comparing them to a natural Greedy algorithm. 
\end{itemize}

The remainder of this work is organized as follows:
In chapter \ref{literature} we shortly summarize the related literature. 
In chapter \ref{formulation} we introduce the problem formally.
In chapter \ref{localsearch} we give a short introduction into different local search algorithms.
In chapter \ref{heuristic} we present a heuristic for the VRD using nested local searches. 
In chapter \ref{computation} we discuss computational results and the impacts of different parameters.
In chapter \ref{outlook} we give an outlook on possible variants and future aspects of research. 
In the appendix we provide a proof for theorem \ref{almostfeasiblethm} stated in chapter \ref{formulation}. Also we provide detailed information about the computational results.

\section{Literature Discussion}
\label{literature}

The problem we are dealing with is related to the field of vehicle routing problems. The common idea of all these problems is that there is a fleet of trucks given which has to deliver a given set of packages to certain positions. Mostly there are also some kind of constraints for the trucks, e.g. they are not allowed to deliver more than a fixed number of packages. 
In this investigation we do not consider packing restrictions yet, which would be an interesting generalization and should be considered in future research. For good books on the vehicle routing problems see e.g. \cite{vehiclerouting} and \cite{golden2008vehicle}. For more information on the vehicle routing problem see also \cite{vehiclerouting2}, \cite{vehiclerouting4}, \cite{vehiclerouting5}, \cite{vehiclerouting6}. However, the problem we are dealing with differs not only superficially from "standard" vehicle routing problems, as in our model synchronization of drone and truck tours becomes a crucial aspect. A step in this direction was already done in \cite{drexl1}, however, the "trailers" which correspond to drones in our setting cannot move by themselves. The consideration of drones within combinatorial optimization problems has started very recently. Murray and Chu introduced such a model of the TSP with drone in \cite{murray}.
They call the problem the Flying Sidekick Travelling Salesman Problem (FSTSP). Other work related to the topic is \cite{agatz}, \cite{ha2}, \cite{ponza},  \cite{ferrandez}, \cite{wang}, \cite{chen}, \cite{dorling}. We (very roughly) summarize three of these papers in the remaining chapter.

\subsection{``The flying sidekick travelling salesman problem: Optimization of drone-assisted parcel delivery'' by Murray and Chu (\cite{murray})} 
In the paper by Murray and Chu \cite{murray} the synergy of drones and one truck was investigated.
Two different models are introduced; for both models, heuristics and mixed integer linear programs are provided. The Flying Sidekick Travelling Salesman Problem is one of these two models, where one truck and one drone (which has restricted flight endurance) are supposed to deliver packages. There are also some subtleties, e.g. there are some packages which cannot be carried by a drone (e.g. because they are too heavy). Note that in our model we abstract from such subtleties.
The objective is to reduce the latest return of the truck or drone after the delivery of all packages.
They provide an integer linear programming formulation which they try to solve for modest-sized instances (10 packages to deliver) with the solver "Gurobi". However, even for these comparatively small instances, "Gurobi" did not find provably optimal solutions within 30 minutes. Thus heuristics seem necessary. They provide a heuristic which basically starts with a solution where the truck delivers all packages on its own. Then, very roughly, in each iteration several changes (using also the drones) are investigated and the change with the highest saving is applied. Computational results for the heuristic are provided.
The second model they propose contains several drones (and still one truck). However, the drones do not interact with the truck any more, but can only start from the depot, deliver a package within their reach and then return to the depot. Here, unlike in the previous model, drones may leave the depot several times. The objective remains the same. Thus here one has to find a truck tour that covers at least the packages which cannot be serviced by the drones (e.g. because of their weight or because they are too far away from the depot). The packages which are not delivered by the truck have to be serviced by the drones from the depot directly. So therefore a  schedule for the drones has to be found.
Here also a heuristic is proposed. Very roughly, it starts with a solution in which all packages which can be delivered by drones in fact are delivered by drones; the remaining packages are delivered by trucks. Then a local search heuristic is applied to improve this solution. 
In both models the computational analysis emphasized the importance of using a good algorithm to solve the travelling salesman problem, which is a subroutine in the proposed heuristics.
Considering multiple trucks and drones is suggested as a topic for future research (which is what we do in this paper).

\subsection{``Optimization approaches for the travelling salesman problem with drone'' by Agatz, Boumann and Schmidt (\cite{agatz})}
The problem which is investigated in \cite{agatz} is very similar to the already discussed problem from \cite{murray}. The main difference is that here "the truck and the drone travel on the same road network" \cite{agatz} (as opposed to the case when the drone is moving according to the Euclidean metric while the truck is moving according to the Manhatten metric like in \cite{murray}).

After introducing the problem, some theoretical aspects are considered, e.g. it is shown that if the truck has speed $1$ and there is one drone with speed $\alpha$, then the optimal solution to the truck only version (i.e. the TSP) is at most $1+\alpha$ times the optimal objective value of the truck drone combination version.

Also an example where this approximation factor occurs is provided, which shows that the factor of $1+\alpha$ is tight. Another application from that fact is that the TSP-D (as the problem is called there) is constant-factor approximable in polynomial time (take e.g. the Christofides heuristic and obtain a $1.5+1.5\alpha$ approximation factor) and this approximation factor is further improved. Afterwards a mixed integer programming formulation is provided. Then heuristic approaches are developed, based on so called route first - cluster second procedures (see \cite{routefirst}).

Then a computational study on artificially generated instances is performed.
Moreover it is suggested to investigate the setting with multiple trucks and drones.

\subsection{``The vehicle routing problem with drones: several worst-case results.'' by Wang, Poikonen and Golden (\cite{wang})}
In \cite{wang} a vehicle routing problem with drones is investigated and several worst-case bounds are proved.
One central proof idea for obtaining worst case results is that given a solution for trucks and drones, a single truck can drive along all tours of the trucks and drones in the given solution one after the other (the resulting tour is then canonically shortened) obtaining a solution using only one truck.
 Like in our model, a fleet of trucks is considered where each truck is able to carry several drones. One difference between their and our model is that in their model a drone never changes the truck it belongs to which is possible in our model.

\section{Problem Formulation}
\label{formulation}

Now we want to introduce the problem formally.
\begin{table}
\caption{}
\label{notation}
 \begin{tcolorbox}[tab2,tabularx={p{1.3cm} p{2cm} p{4cm}},title=Data defining an instance]

  \textbf{symbol} & \textbf{domain} & \textbf{meaning} \\ \hline
  $n_t$ & $\N$ & number of trucks  \\ \hline
  $n_d$ & $\N$ & number of drones  \\ \hline
  $n_p$ & $\N$ & number of packages \\ \hline
  $Pos$ & $(\Z \backslash \{(0,0)\})^{n_p \times 2}$ & positions of packages \\ \hline
\end{tcolorbox}
\end{table}
\subsection{Instance}
Any instance consists of the data described in table \ref{notation}, which is the number of trucks, the number of drones, the number of packages and the destination positions of the packages. 
Every truck can carry an arbitrary number of packages and drones and the depot is always $(0,0) \in \Z^2$. 
The packages together with the depot are the nodes of a complete graph, e.g. from each package (or the depot) one can travel to any other package (or the depot). Whenever we refer to subgraphs we refer to subgraphs of that complete graph.
This can also be the case implicitly, e.g. if we refer to a "cycle which contains the depot", then we refer to such a cycle which is a subgraph of the introduced complete graph.
While this complete graph is undirected, we also consider directed subgraphs in the sense that the edges of the subgraph are a subset of the edges of the complete graph, but they have assigned an orientation.
Later we will also define certain multi graphs on the same node set. We simply enumerate the drones by $1,...,n_d$ respectively and analogously for trucks and packages; e.g. we write either $d_i, t_j, p_k$ or $i,j,k$ to refer to certain drones, trucks or packages respectively, but we only use the second notation if the meaning is clear from the context. 

\subsection{Solution} 
We now define a solution formally. Note that not every solution will also be a feasible solution. We will define feasibility afterwards. First we need some preliminary definitions.
\begin{definition}[Truck  / Drone Tour]
A tour is a directed cycle $C$ (in particular all nodes are distinct) which contains the depot and at least one node from $\{1,...,n_p\}$.

A truck tour is a tour with a function "carry" defined on the edges $E$ of the truck tour, $\text{carry}: E \rightarrow Pot(\{1,...,n_d\})$, which indicates the drones which are carried by the truck on the respective edge. Similarly a drone tour is a tour with a function "carried" defined on the edges $E$ of the tour,
$\text{carried}: E \rightarrow \{0,1,...,n_c\}$, indicating on which truck the drone rides (0 means it flies on its own). We call all functions "carry" and "carried" the "carry functions".
\end{definition}

Each package is delivered by the first vehicle arriving at the node which is not a drone that already has delivered a package on its current flight.

\begin{definition}[Solution]
A solution consists of a truck tour for each truck ($T_{t_1},...,T_{t_{n_t}}$) and a drone tour for each drone ($T_{d_1},...,T_{d_{n_d}}$). 
\end{definition}

\subsection{Feasibility of a Solution}
If we just look at any solution it is not necessarily practicable. A solution is feasible if the following points hold:
\begin{itemize}
\item  Each node (except from the depot) is visited at most either by a single drone or by a truck as well as the drones it carries when arriving and the drones it carries when leaving. 
\item Vehicles that arrive at the depot stay there.
\item Each node is visited by at least one vehicle (truck or drone).
\item The number of consecutive zero edges in a drone tour is bounded from above by two (we don't count around the depot here), but the drone has no restriction in flying time or distance. 
\item The solution has to be consistent.
\end{itemize}
\begin{definition}[Consistency]
A solution that fulfils the other conditions of feasibility is consistent if the following two conditions hold.
\begin{itemize} 
\item (Carry Consistency) For each edge $e$ of a truck tour of a truck $t \in \{1,...,n_t\}$ with $\text{carry}(e)=D$, the respective drone tour $T_d$ for each $d \in D$ contains also the edge $e$ and the carry function of the tour $T_d$ fulfils $\text{carried}(e)=t$. Conversely, for each edge $e$ of a drone tour from a drone $d \in \{1,...,n_d\}$ with $\text{carried}(e)=t \neq 0$, the respective truck tour $T_{t}$ contains the edge $e$ and the carry function of the tour $T_{t}$ fulfils $d\in \text{carry}(e)$.
\item (Schedule Consistency) If the above condition is fulfilled we can define a relation on the edges: 
Two edges are essentially equal if the corresponding vehicles travel this edge together. Because of the above condition this leads to a partition of all edges of all tours of a solution. Edges that are mapped to $0$ or $\emptyset$ by the respective carry function are only essentially equal to themselves.
Now the second condition can be formulated as follows:
Look at a solution, i.e. all the drone and truck tours of the solution. 
In the first step we mark all edges leaving the depot.
In each further step we define $E^*$ to contain those edges for which all the previous\footnote{W.r.t. the order in the directed cycle they are part of, starting from the edge leaving the depot.} edges are already marked; in each such step we mark the following edges of $E^*$: Every edge for which all essentially equal edges of it are also in $E^*$. This process has to lead to a completely marked solution (i.e. every edge of every tour is marked).
\end{itemize}
\end{definition}

\begin{figure}
\includegraphics[scale=1]{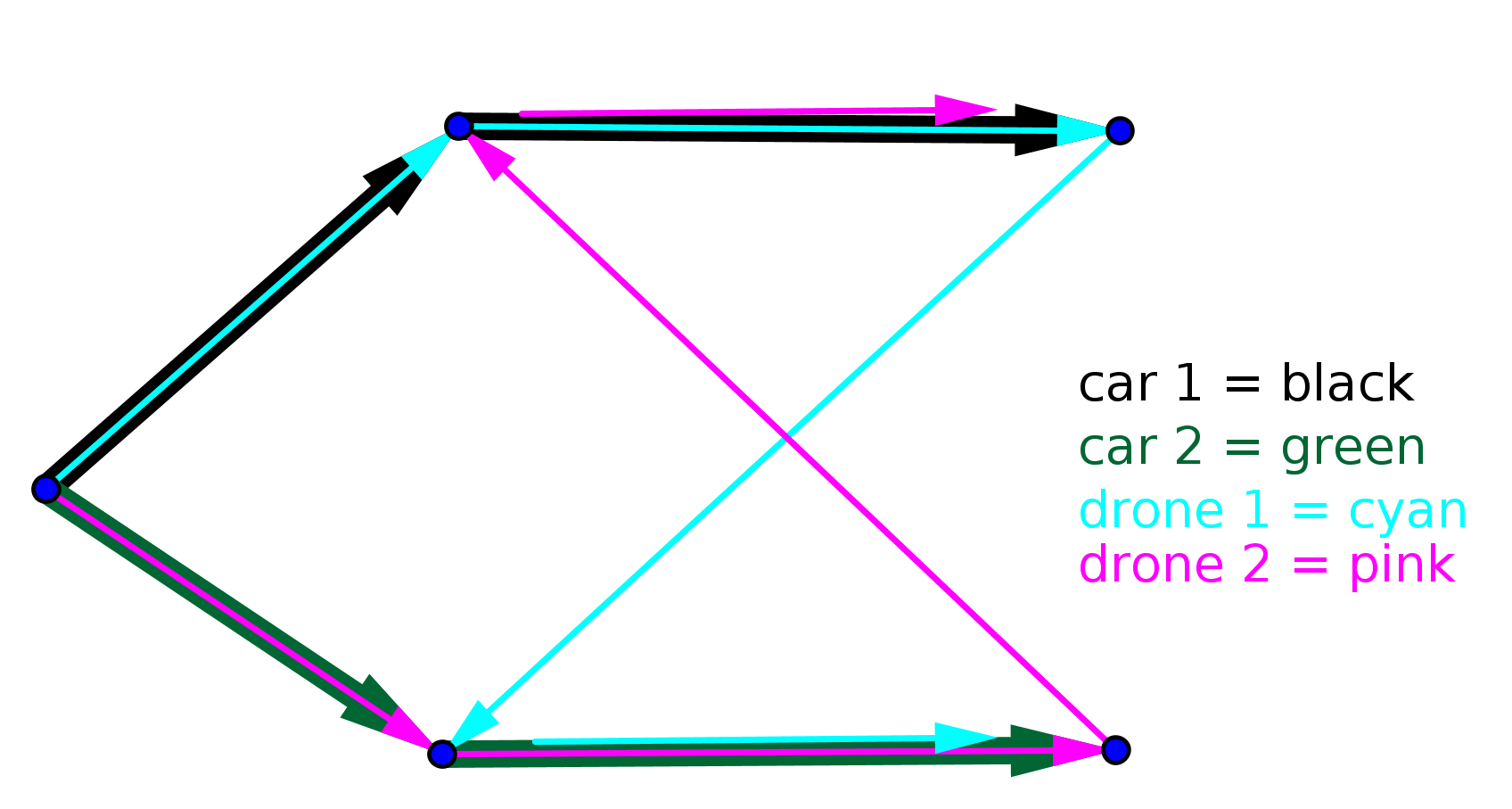}
\caption{An example of inconsistency}
\label{nonconsistency}
\end{figure}

\begin{figure}
\includegraphics[scale=1]{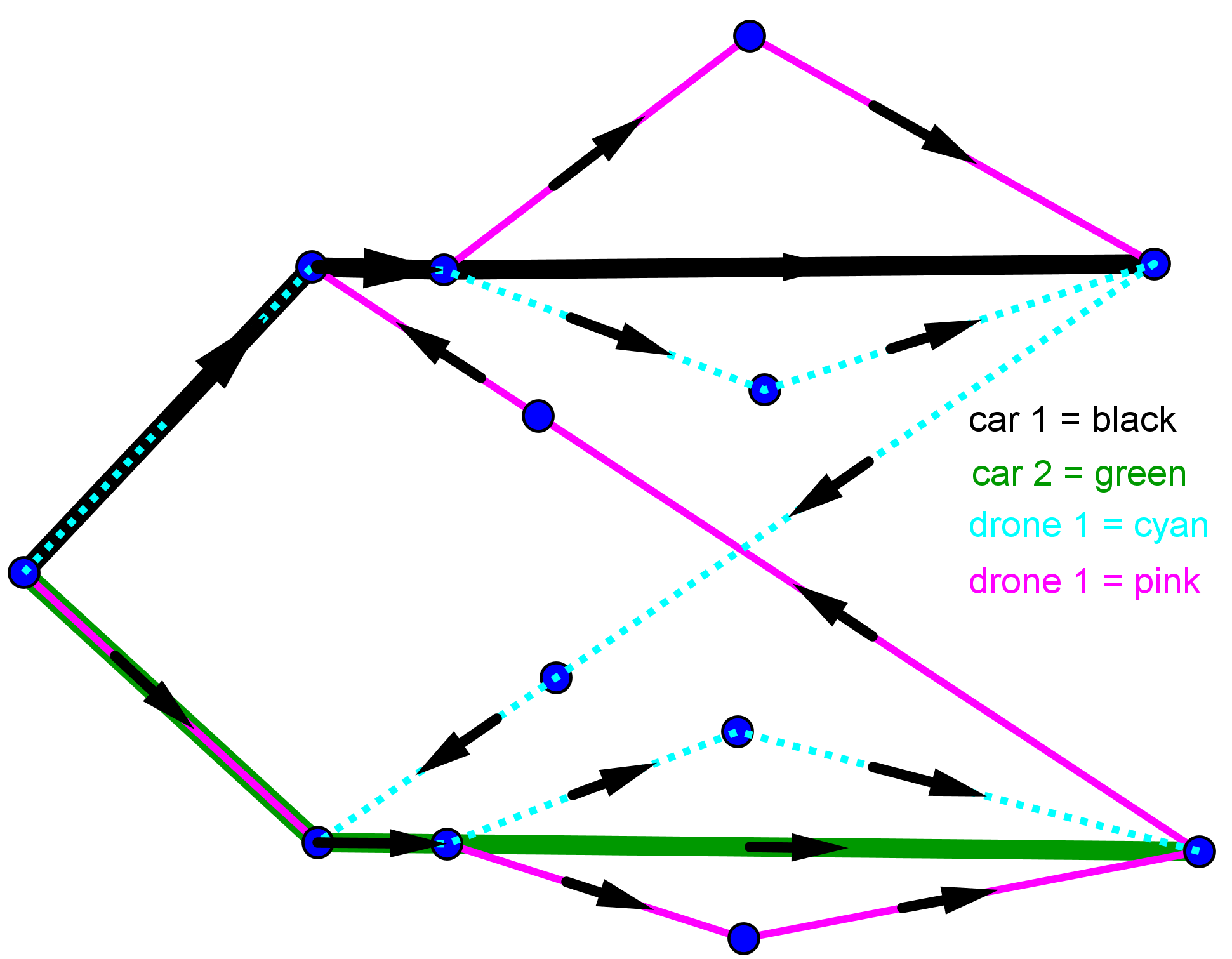}
\caption{A more severe example of inconsistency}
\label{nonconsistency2}
\end{figure}

Examples of inconsistent solutions can be found in \ref{nonconsistency} and \ref{nonconsistency2}.

\subsection{A Characterization of Consistency}
Here we provide a new characterization of consistency.

\begin{definition}[Graph of a Solution]
We obtain the directed multi-graph of a solution as follows:
The nodes are $0,...,n_p$. The edges $E$ are the union of the edges from all truck and drone tours. We also define analogously to above a carry function on the edges of the graph, that is a canonical extension of every carry function on a single tour. Moreover we have a function that returns for every edge to which vehicle it corresponds (e.g. drone number 3).
\end{definition}

\begin{definition}[Almost Feasible]
We call a solution almost feasible if everything needed for feasibility except the second condition of consistency, the schedule consistency, is fulfilled.
\end{definition}

\begin{theorem}
\label{almostfeasiblethm}
Each almost feasible solution is consistent (and thus feasible) if and only if its graph has only cycles which have two consecutive drone-only edges representing different drones or contain the depot. We call such cycles flip-cycles.
\end{theorem}
The proof can be found in the appendix.

\begin{figure*}
\includegraphics[scale=1]{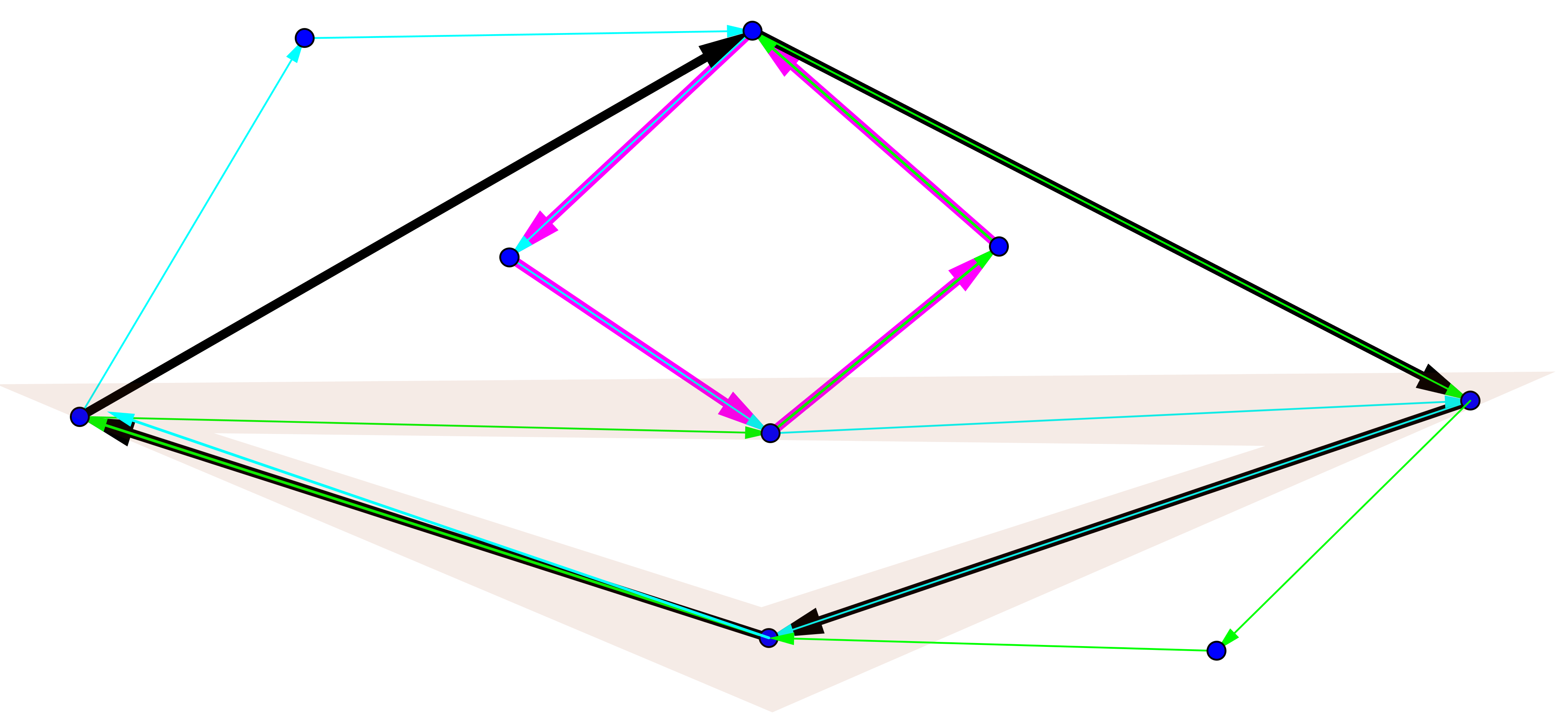}
\caption{It is necessary to allow for flip-cycles.}
\label{strongversionfalse}
\end{figure*}

Note that the theorem would be wrong, when one does not allow for flip-cycles. For a counterexample see figure \ref{strongversionfalse}.

\subsection{Objective}

We are left to define what the objective function is. We assume that the drones and trucks are moving at the same speed (this assumption can easily be relaxed as our algorithm does not exploit it), however drones move according to the Euclidean metric while trucks move according to the Manhattan metric. This reflects the fact that trucks are restricted to the street network.
Now the objective is to minimize the average delivery time of the packages.
Recall that a package is delivered by the first vehicle arriving at the respective node which is not a drone that has already delivered a package on its current flight.
For computing the objective value also note that drones and trucks, depending on the solution, often have to wait for each other. \\ \\ 
We note, that this is not the only interesting objective. The literature often analyses the completion time, in particular either the time the last vehicle returns or the average of the returning times of all vehicles. For parcel services it is crucial to know how many trucks and drones they should purchase for their delivery area. This questions is treated implicitly in some of our computations when the numbers of trucks and drones are varied.

 \section{General Local Search Heuristics}
\label{localsearch}
The principle of local search algorithms is to search within a neighbourhood of a current solution for new solutions, which are --- or at least tend to be --- better than the current solution. For many combinatorial optimization problems local search methods are among the state of the art algorithms. For example the Lin-Kernighan heuristic effectively implemented by Keld Helsgaun is among the state of the art algorithms for the TSP, see \cite{lkh1} and \cite{lkh2}. For a precise treatment of the topic of (stochastic) local search heuristics, we refer to the book written by Holger Hoos and Thomas Stützle \cite{hoos1}. The core of each local search algorithm is the definition of the neighbourhood(s). Vice versa, by defining the neighbourhood(s), a large part of a local search method is already defined. Metaheuristic frameworks exploit this fact by letting the user customize the problem instance, the solutions, neighbourhood(s) and initialization procedure and then providing a framework for various local search methods. We used the local search framework JAMES \cite{james} for implementing our algorithm. Because of the fact that a local search can be customized essentially by defining appropriate neighbourhoods, one can put one's focus on defining good neighbourhoods and then one can apply several local search methods using the same neighbourhood definitions. We now shortly explain some of the main principles of several local search algorithms.

\subsection{Steepest Descent --- Random Descent}
Steepest descent is a very basic local search algorithm. From a current solution, the algorithm investigates the complete neighbourhood and moves on to the best new solution if it is better than the current. If no solution is better than the current, then the heuristic will terminate. This heuristic has the major drawback that it has no mechanism to avoid getting stuck in local optima. This can be mitigated by rerunning the algorithm with different initial solutions; however, this often will not suffice, especially when there are many local minima. Random descent is a very similar heuristic, but instead of searching the whole neighbourhood for the best solution, in each step a random neighbour is chosen and if it is better, it is the new solution, otherwise a new random neighbour is created. One has to define a suitable stopping criterion, e.g. no improvements for a certain number of iterations. Just as steepest descent, this heuristic might get stuck in local optima.
Of course there is a canonical mixture from random descent and steepest descent: Just create $k$ random neighbours and take the best if it is better than the current solution.

\subsection{Tabu Search}
Tabu search can be thought of as an attempt of taking the positive aspects of steepest descent while avoiding its main disadvantage, i.e. the disability of escaping from local optima. Tabu search was introduced by Glover (compare \cite{tabu0}, \cite{tabu1} and \cite{tabu2}). A beautiful book by Glover and Laguna gives a good overview and introduction into tabu search \cite{tabu3}. Moreover the tabu search tutorial from Glover from 1990 \cite{tabu4} as well as Laguna's guide \cite{tabu5} are recommendable. Here we are only going to provide a rough sketch of tabu search, omitting many details and variants.

\subsubsection{Short-Term Memory for Escaping Local Extrema}
A very basic tabu search is a simple generalization of the steepest descent algorithm: Basically one performs the steepest descent algorithm and always remembers the best solution so far. However, here the current solution is not necessarily the best solution. In the end, the best solution found is returned. The current solution is in principle computed like in the steepest descent heuristic. However, to be able to escape from local minima one keeps track of a number (can be always the same fixed number, but can also vary over time) of solutions that one has visited last. It is forbidden to go there again, even if it is the best solution in the neighbourhood of the current solution. Then one goes to the best solution from the neighbourhood where one is allowed to go. This enables the search to escape from local optima. Note that tabu search in its general form is not restricted on forbidding certain solutions. More general, moves can be forbidden, which for example for the TSP problem could be to forbid swapping cities $i$ and $j$. In these cases so called aspiration criteria can improve the results. If such an aspiration criterion is fulfilled (e.g. the respective move yields a solution that is better than the best solution found so far) it is accepted as the new current solution even if it is tabu. The list where the algorithm keeps track of the previous solutions or moves respectively is often referred to as the short-term memory.

\subsubsection{Intermediate-Term Memory for Intensification}
 The intermediate-term memory tries to push the current solution into good regions. Therefore --- loosely spoken --- it looks at good solutions that are collected so far and extracts some features of them. E.g. in good TSP solutions, only a small subset of the edges of the graph might be used. Thus solutions using these edges might be rewarded and thus might even be preferred over solutions which have a better objective value. Compare again \cite{tabu1}. Thus this procedure causes an intensification of the search.
 
 \subsubsection{Long-Term Memory for Diversification}
 The long-term memory has a complementary function compared to the intermediate-term memory. Its purpose is to diversify the solution by penalizing features that occurred often. For example for edges that occurred in many TSP solutions found so far (not only in the good ones like above) there may be a penalty for solutions containing them. This diversifies the search, as new parts of the search space are visited.

\subsection{The Metropolis Algorithm}
Basic Metropolis search is a heuristic based on the more general break-through works \cite{metro1} and \cite{metro2}. Here we have a fixed temperature and a solution. Then we look at a random neighbour and accept it if it is better or accept it with a certain probability if it is worse. This probability depends on how much worse it is (much worse solutions are unlikely to be accepted) and on how high the temperature is (higher temperatures correspond to higher probabilities). Metropolis search can be considered as a simplification of simulated annealing, but there, the temperature is decreasing; compare \cite{simannealing1}, \cite{simannealing2} and \cite{simannealing3}.) 
or as Ingo Wegener states in \cite{simannealing4}, the Metropolis algorithm is equivalent to simulated annealing without temperature changes.

\subsection{Replica Exchange Monte Carlo --- Parallel Tempering for Combinatorial Optimization}
This method was introduced by Swendsen and Wang \cite{tempering1} and is based on Metropolis Search. See also \cite{tempering2}, \cite{tempering3} and \cite{tempering4}.
However the previous mentioned literature does cover parallel tempering with no special focus on combinatorial optimization. For its application to combinatorial optimization we refer to \cite{tempering5} where the application of parallel tempering to the TSP is illustrated. We roughly summarize the principle of parallel tempering for combinatorial optimization staying close to \cite{tempering5}:
The core idea is to look at several copies of a system, assigning each a temperature and then changing the system depending on the temperature (where higher temperatures make bigger changes more probable, compare again the Metroplois algorithm). However, the systems are not considered separated: As already said, parallel tempering simulates multiple parallel solutions, each such system having assigned a certain temperature. Let us say there are $k$ such solutions, with temperatures $T_1<...<T_k$. Each solution is modified according to the basic Metropolis search described above. After some iterations, randomly two consecutive (with respect to the temperature) solutions are swapped (i.e. their temperatures are exchanged) with a certain probability.  For details we again refer to \cite{tempering5}.  Because of this principle, solutions can raise in temperature thus overcoming local minima and then, when cooled down again, the solutions can step into new, hopefully better, local minima.

\subsection{Piped Local Search}
Piped local search is a very simple and natural idea mentioned in \cite{james}. One simply applies several local search heuristics one after the other, each time taking the solution of the previous local search as the new initial solution.

\subsection{Summary}
The heuristics are summarized in table \ref{heuristics}.
\begin{table}
\caption{}
\label{heuristics}
 \begin{tcolorbox}[tab2,tabularx={p{3cm}p{4.7cm}},title=Local Search Heuristics]
  \textbf{Heuristic} & \textbf{Key Idea}  \\ \hline
  Steepest Descent & Define neighbourhood. Take initial solution. Look at entire neighbourhood. If one solution in this neighbourhood is better, replace the current solution with the best solution of the neighbourhood, otherwise terminate. Iterate this.    \\ \hline  
  Random Descent & Basically like steepest descent, but instead of looking at all neighbours, one only looks at one randomly sampled neighbour and accepts it if it is better. This process is iterated. Clearly there is a canonical variant in the middle of random descent and steepest descent. \\ \hline
  Tabu Search & Tabu search can be considered as a generalization of the steepest descent heuristic. However, it has a mechanism allowing to escape from local minima, i.e. a tabu list of, e.g., moves depending on moves made recently (short-term memory). It also has mechanisms for diversification and intensification.  \\ \hline 
  Metropolis Search & Fix a temperature. Choose a solution. Look at a neighbour. Accept if it is better. Accept also with a certain probability if it is worse, depending on how much worse it is and depending on the temperature.\\ \hline 
  Parallel tempering & Looks at several copies of a system, each having a different temperature. Performs Metropolis search on each copy. After some time, e.g. periodically, two neighboured systems (with respect to temperature) switch their temperatures. \\ \hline 
  Piped Local Search & Simply apply different local search methods in a row, each time taking the previous solution as the new initial solution.   
\end{tcolorbox}
\end{table}

\section{The Algorithm} \label{heuristic} 
\begin{figure*}
\includegraphics[scale=0.9]{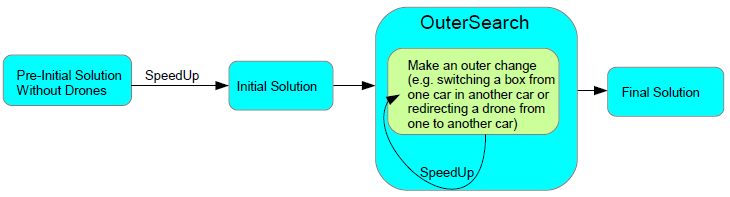}
\caption{Basic structure of the algorithm}
\label{basicstructure}
\end{figure*}

Now we describe the algorithm we developed, which consists of two nested local search procedures:

\begin{enumerate}
\item As a pre-initial solution, solve the multiple TSP, i.e. solve the problem only using trucks and ignoring the drones.
\item Assign an appropriate number of drones to each truck, which, in this step, stay associated to this truck. Then, for each drone, use a local search approach (called SpeedUp) to disburden the truck it is associated to as much as possible. However, note that SpeedUp works locally, in particular, packages are not switched between trucks and drones are not redirected to other trucks. This is the initial solution. 
\item After that, an outer local search (called OuterSearch) starts, which does not treat the tours isolated. In contrast to SpeedUp, packages can be interchanged between the trucks to obtain new solutions and moreover drones can switch between trucks. Note that SpeedUp can also be applied to a part of a truck tour instead of applying it on the whole tour, which is done after a drone switch. 
\end{enumerate}
The basic structure of the algorithm is illustrated in figure \ref{basicstructure}.

Now we describe the steps of the algorithm in more detail.

\subsection{Initial Solution}

Given an instance as defined earlier, we want to find a feasible solution which is supposed to be our initial solution.
Therefore we apply the following steps.

First, we solve the multiple TSP (mTSP) with our objective function, i.e. the TSP with multiple salesman. This is our pre-initial solution, in particular until now all drones are useless. To solve the mTSP we 
 assign the packages sorted by their angle (around the depot in $(0,0)$) to the trucks. Then we solve the TSP (with a simple local search provided in the illustrating examples of JAMES \cite{james}, we just changed the objective function) for each truck. This is repeated for several start angles, the best result is used as pre-initial solution.

Now we equally distribute the drones to the truck tours from the pre-initial solution. Each drone first rides on its truck the whole time, and then a local search algorithm called SpeedUp (defined in \ref{speedup}) is applied for each drone. The result is our initial solution.

\subsection{SpeedUp}
\label{speedup}
\begin{figure}[H]
\includegraphics[scale=0.3]{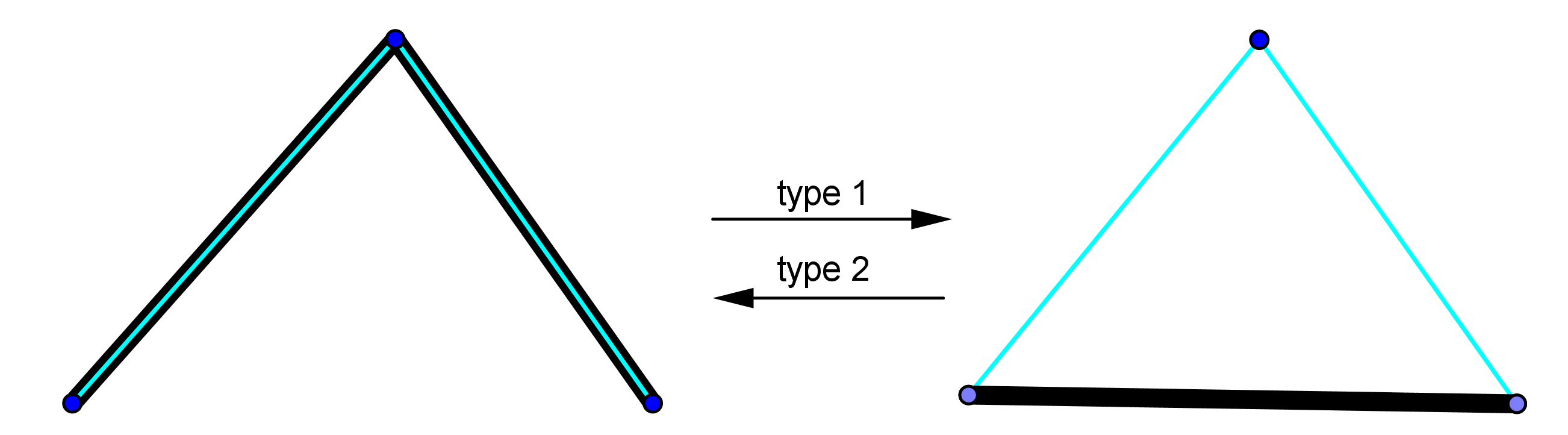}
\caption{Small Neighbours of type 1 and 2}
\label{n1}
\end{figure}

\begin{figure}[H]
\includegraphics[scale=0.5]{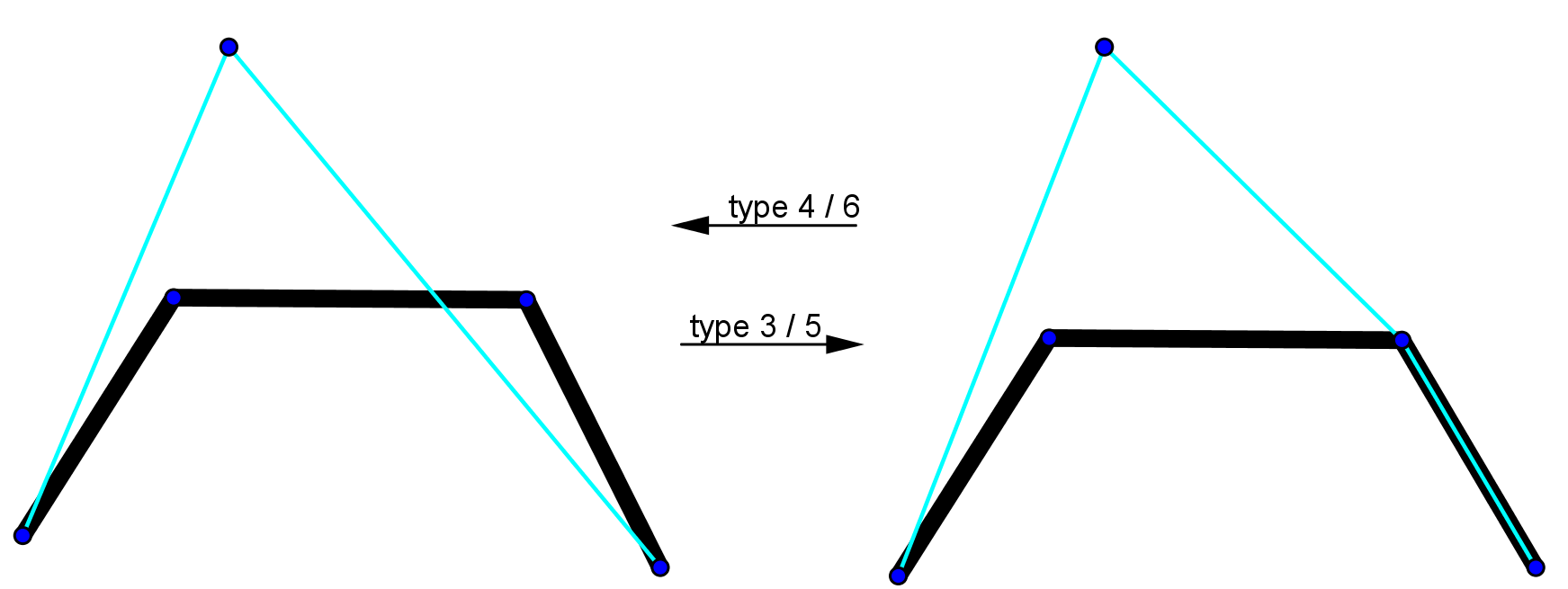}
\caption{Small Neighbours of type 3,4,5 and 6}
\label{n2}
\end{figure}

SpeedUp acts on (a part of) a drone tour and the associated truck tour that is interwoven with that drone (on this part of the drone tour). The drone may not change the truck (in this part of the drone tour). The part of the overall solution, where this drone and its assigned truck interact, is called SpeedUp area. 
SpeedUp is a local search approach that aims for improving the tour by using the drones as efficiently as possible.
As the core of a local search is the neighbourhood definition, most importantly we have to explain the neighbourhood definition of SpeedUp.
Some key factors of a solution are not changed by SpeedUp, i.e. the assignment of drones and packages to trucks stays the same in the SpeedUp area. 
We define a neighbour as a solution that arises by applying one of the following nine kinds of moves. 
While applying the nine types of changes, two kinds of problems may appear:
\begin{itemize}
\item[$P_D$] The drone is unable to fly that far.\footnote{Note that in our computations we will consider the case where the drone has an unlimited range. However, our algorithm does not depend on that and can easily be used in the case with a restricted drone range.}
\item[$P_T$] The truck misses a node that he has to visit because he either receives or sends away another drone there.
\end{itemize}
If one of the above problems occurs, this change is simply omitted and then a new change is tried. 
We now look at these possible changes in detail. As opposed to the neighbours that we are going to define in OuterSearch later, these neighbours in the local search SpeedUp performs are called small neighbours.
\paragraph{Small Neighbours of Type 1 and 2}
Small neighbours of type 1 arise through the following change of a tour of a solution: If there are two edges on which a truck drives carrying a drone, this can be changed by letting the truck drive from the start node from the first edge to the end node of the second edge. The drone delivers the remaining package in the middle. Here both problems, $P_D$ and $P_T$, can arise. Compare figure \ref{n1}.\footnote{In all figures illustrating the different moves, the edges are directed from left to right if not stated differently.} Note that in all figures, the drone's path is indicated in cyan while the truck's path is indicated in black. Small neighbours of type 2 arise by using the inverse move, without having to care about $P_D$ or $P_T$. 
\paragraph{Small Neighbours of Type 3,4,5 and 6}
Small neighbours of type 3 are obtained by changes of the following kind: 
If the drone leaves the truck, delivers a package and then returns to the truck later, in the neighboured solution it returns earlier to the truck and then stays there at least until arriving at the node where it would have returned usually.
Here only problem $P_D$ may arise. Small neighbours of type 4 arise by using the inverse move, here also only $P_D$ may arise. Compare figure \ref{n2}.
 Small neighbours of type 3 and 4 describe changes when the return to the truck is performed earlier or later. Conversely we obtain neighbours of type 5 and 6 by performing the departure earlier or later. Figure \ref{n2} also illustrate type 5 and 6 if one simply imagines that now the edges are travelled in the opposite direction as one imagined before. Again only $P_D$ may occur.
 
\begin{figure}[H]
\includegraphics[scale=0.65]{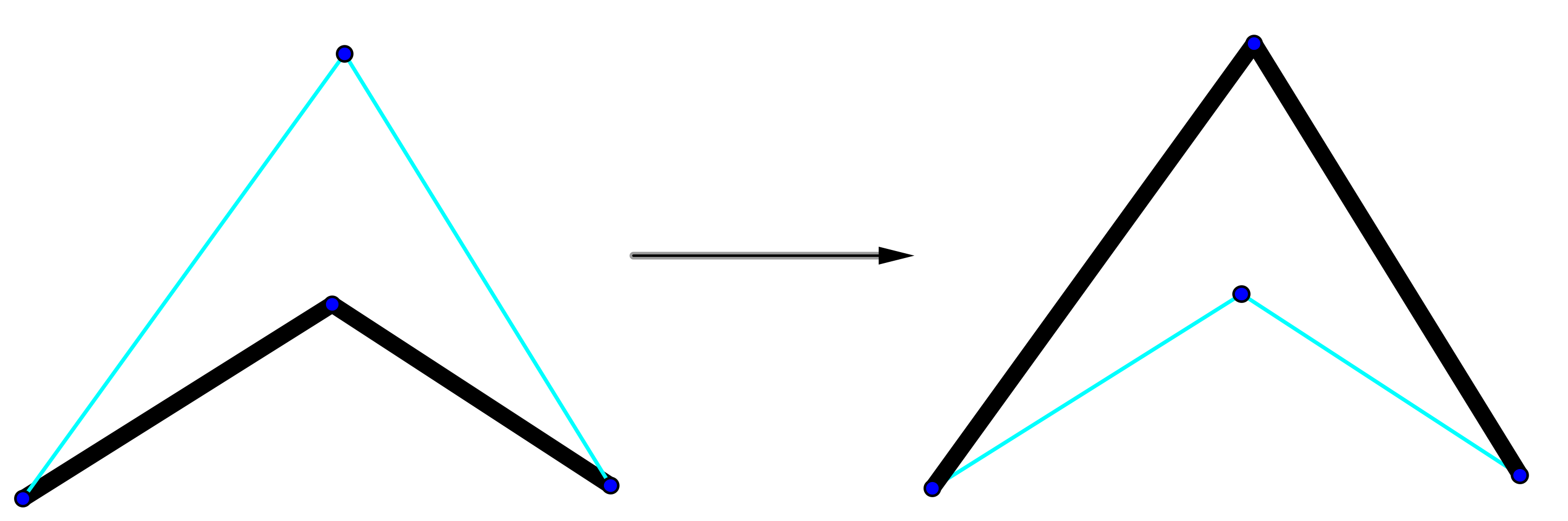}
\caption{Small Neighbours of type 7}
\label{n4}
\end{figure}

\paragraph{Small Neighbours of Type 7}
The drone leaves the truck, delivers a package and returns to the truck; meanwhile the truck has delivered exactly one package in between. Now we obtain a new neighbour by flipping the truck's and the drone's ways. Here problems $P_D$ and $P_T$ can arise. Compare figure \ref{n4}.

\begin{figure}[H]
\includegraphics[scale=0.6]{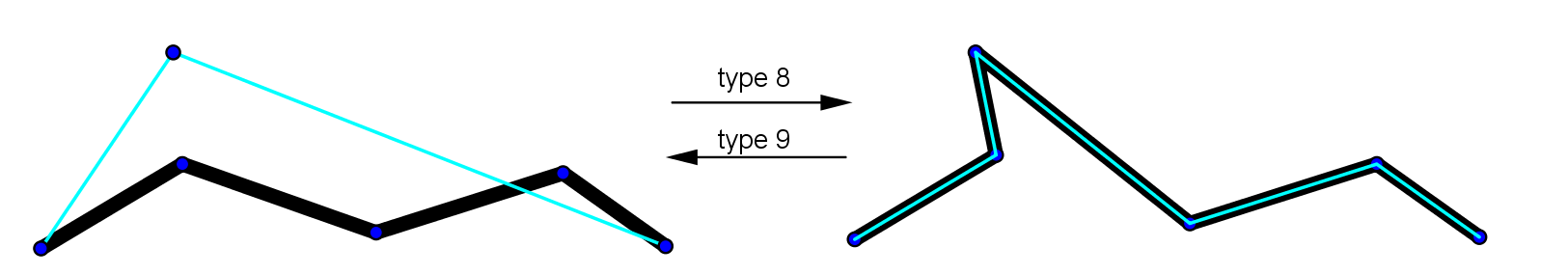}
\caption{Small Neighbours of type 8 and 9}
\label{n6}
\end{figure}

\paragraph{Small Neighbours of Type 8 and 9}
Neighbours of type 8 are especially simple, as here neither $P_D$ nor $P_T$ can occur.
Here the initial situation is that a drone is starting from a truck, delivering a package and then returning to the truck. Now in the neighboured solution the drone does not leave the truck, but the truck delivers the package instead, either before or after the truck's nearest (resp. Manhattan metric) position, whatever of both takes less time for the truck. 
Compare figure \ref{n6}. Note that here only an immediate improvement could occur if the trucks were faster than the drones. However, even when there is a worse solution, this can be useful to overcome local minima.
We obtain neighbours of type 9 by performing the inverse, which is a bit more tricky. If the drone is carried by a truck over enough edges, the package which is left to the drone is chosen uniformly under the possible ones considering $P_T$. Afterwards the drone leaving node and the drone arrival node are also chosen uniformly under the possible ones considering $P_D$. 
SpeedUp now works, as previously stated, on a part of a drone tour, where the drone is just interacting with one truck, called SpeedUp area. A random neighbour is chosen with the following procedure: 
\begin{itemize} 
\item Choose a random \footnote{We always use uniform distribution for randomness, more advanced random distributions could be an aspect of future research.} node of the drone tour in the SpeedUp area. 
\item Check which small neighbours are possible there. 
\item Choose a random small neighbour which is possible and apply it. 
\end{itemize}
\subsection{OuterSearch}
Now, with the help of the SpeedUp function, we can define a procedure that tries to find a new neighbour in the OuterSearch. 

It consists essentially of the following steps (given a current solution): 

\begin{enumerate}
\item Flip a fair coin to decide whether to apply a drone change or a package change. 
Investigating more elaborate decision mechanisms could be an aspect of future research.
\item If a drone change was selected, do the following:
\begin{enumerate} 
\item Choose randomly a drone, a truck, a node in the drone tour (leaving node) and a node in the truck tour (arriving node). Note, that it's possible, that the drone directly returns to the depot (which would be the last node of every truck tour) or that it starts with the truck right from the beginning (as the leaving node and the arriving node could be the depot, where drone and truck start). 
\item Distribute the drone's jobs after the drone leaving node: The trucks to which the drone was assigned now service the packages which the drone originally was supposed to service. In particular, if the drone would have left the truck at some node to deliver a package, the truck now delivers this package right after visiting this node and continues its planned schedule afterwards. 
\item The drone tour stays the same until the leaving node, unless the last edge can be skipped, meaning the drone was flying and doesn't have to deliver a package at the end. If this is the case, delete this edge. 
\item Insert the drone at the chosen position of the new truck tour and adapt all involved tours canonically, i.e. the drone travels along with the new truck all the time without delivering any packages. In particular it never leaves the truck. 
\item Perform SpeedUp on the concerning part of the new truck tour, that is after the arriving node of the drone. Therefore, the drone helps the new truck to deliver its packages. 
\end{enumerate}

\item If a package change was selected, do the following: 

\begin{enumerate}
\item Choose a random package $p$ that is delivered by a truck called "fromTruck" (not by a drone) and choose a random truck called "toTruck" that shall deliver this package now. Note, that this could be the same truck, which can result in a new route of the truck. 

\item \textbf{Remove $p$ from the tour of fromTruck}: Here we have to distinguish between two cases.
\begin{enumerate}  
\item[\textbf{Case 1}] \textbf{Package not at drone dropping/taking position} \\
$p$ is not at a position where the truck drops or takes a drone. In this case, we delete the package and fromTruck just skips travelling there. 
\item[\textbf{Case 2}] \textbf{Package at drone dropping/taking position} \\
In this case we cannot simply delete $p$ from the tour. As we do not want to travel there without actually delivering a package, we try a work around. If it fails, we discard this approach of finding an adjacent solution.\\
Work around:
The drop / arrival is instead performed at an adjacent node, if feasible. If the drone reach is not limited, the arising problems still can be that the drone misses a charging edge and arrives at a node, where it should leave immediately. 

If the work around yields a feasible solution, it is performed and the tours of fromTruck and the drones are updated. 
\end{enumerate}
 
\item \textbf{Insert the package in the tour of toTruck}: 
The package $p_{\text{close}}$ from the tour of toTruck which is closest to the destination of $p$ is considered. Then $p$ is inserted either before or after $p_{\text{close}}$ in the tour of toTruck, depending on what saves toTruck the most time.

\end{enumerate}
\end{enumerate}

\subsection{"Greedy Drones"} 
For comparison we use the following natural Greedy algorithm which we call "Greedy Drones" or simply "Greedy". The reason why we call it "Greedy Drones" is that it computes the solution for the mTSP rather extensively exactly like our algorithm, but then, based on the mTSP solution, assigns the drones greedily. 
In detail, it works as follows: 
\begin{enumerate}
\item Compute truck tours that ignore all drones using the mTSP as in the pre-initial solution\footnote{In most cases, we directly took the truck tours from the pre-initial solution which increases comparability; however, for instances with more drones than trucks the pre-initial truck tours had to be recomputed as in the first implementation we used a Greedy approach which was to weak; of course we used the same algorithm and parameters as before to solve the mTSP.} 
and equally distribute the drones to them. 
\item We change each truck tour by now also using the drones. Assume there are $k$ drones assigned to a truck. Then we change the respective truck tour as follows: If at the depot it is possible to send away all drones to service the next (with respect to the order within the original truck tour) $k$ packages, then this is done and the truck meets with all the $k$ drones again at package position $k+1$. Then at least one edge (to package $k+2$) is travelled together, as the drones need to charge, and then the same procedure is repeated. However, if at a node not all drones can fly away, then they travel on the truck to the next node and the respective procedure starts form there; i.e. if it is possible to depart from node $i$ then the $k$ drones depart from there and meet again $i + k+1$ nodes later and so on. Note that in our computations we set the maximal flying distance to $\infty$, thus the drones can always depart if they are charged.

\end{enumerate}

\section{Computational Results and Evaluation} \label{computation} 
 
To investigate our algorithm, we implemented\footnote{We ran our computations on a Microsoft Windows 10 Pro machine with an Intel Celeron CPU G3900 (2.8 GHz, 2 cores) and 4 GB RAM and no dedicated graphic card; the Java version we used is Java 8 Update 111 (Oracle); we used Eclipse (Neon) as our IDE.} it in Java using JAMES. Trying several local search heuristics, using parallel tempering we obtained the best results.   
Thus here we focus on parallel tempering as a local search method. As instances we choose the package positions uniformly at random on a $200 \times 200$ integer grid excluding $(0,0)$.
Note that we used long running times for computing one solution; a more precise description of the setting we used for our algorithm is in the appendix; as we also used a maximum number of steps without improvement as a stop criterion, the running time may sometimes differ for the same settings, details are in the appendix.

We are going to look at several instances with the parameters presented in table \ref{parameters} to obtain a rough overview. 

\begin{table}
\caption{}
\label{parameters}
 \begin{tcolorbox}[tab2,tabularx={p{0.5cm}p{2cm} p{2cm} p{2cm}},title=Different Parameter Settings]
  \textbf{\#} & \textbf{\# Packages} & \textbf{\# Drones} & \textbf{\# Trucks} \\ \hline
  1& 200 &2 &1 \\ \hline  
  2& 200& 2&2 \\ \hline 
  3& 200& 5&3 
\end{tcolorbox}
\end{table}

\begin{figure*}
\includegraphics[scale=0.7]{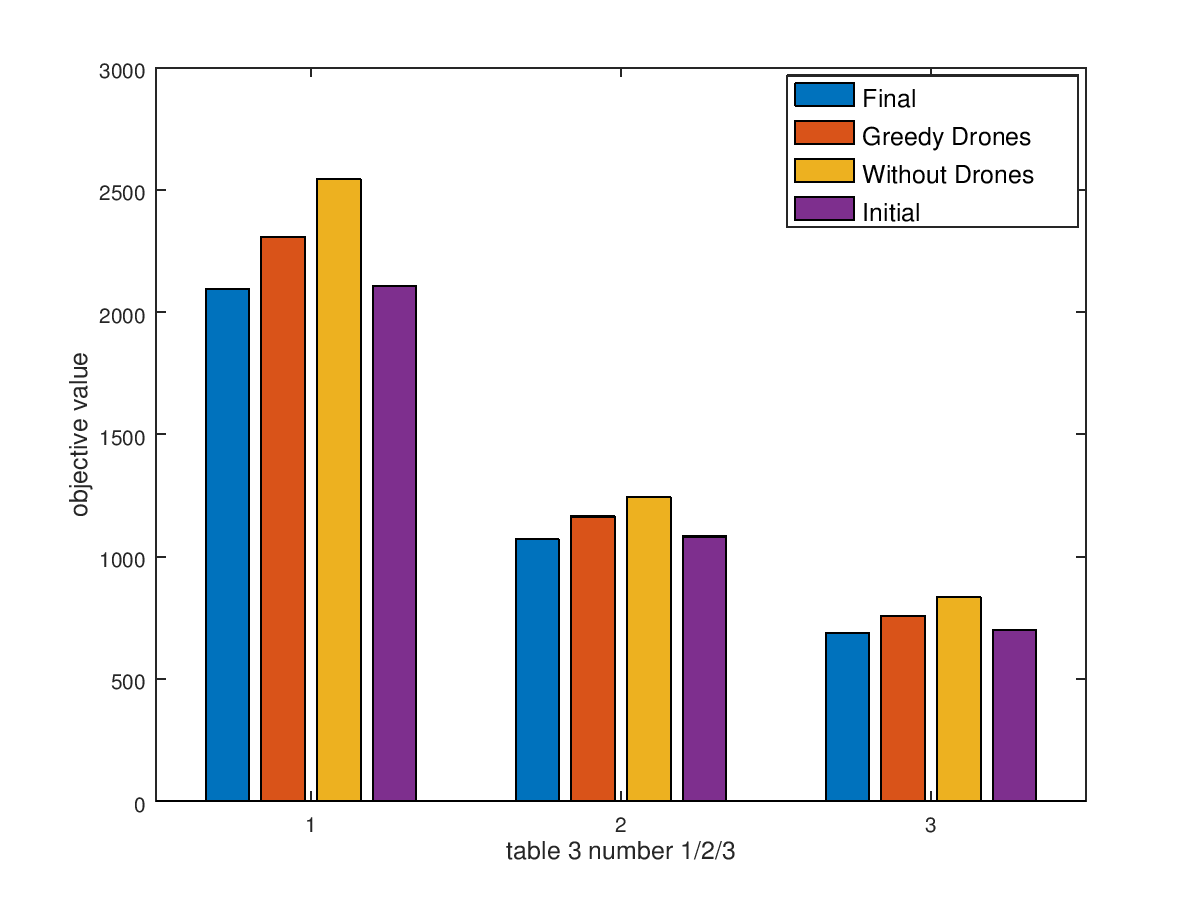}
\caption{Computational Results according to the data from table \ref{parameters}.}
\label{tab3}
\end{figure*}

The results are visualised in figure \ref{tab3}; every setting was averaged over 10 different instances.

So we can see that the solution produced by Greedy Drones yielded results with objective value  10.1\%, 8.7\% and 10.2\%  bigger than the results obtained with our local search approach for table \ref{parameters} numbers 1, 2 and 3 respectively; the solution without drones yielded results with objective value 21.3\%, 16.1\% and 21.0\% bigger than the results obtained with our local search approach, again for table 3 number 1,2 and 3 respectively. We also see that the outer local search has only marginal effect. A possible reason is that we solve the mTSP such that the resulting tours are rather separated (every truck stays in its sector) and that possibly too many of the moves in the outer local search are unlikely to bring improvements to the solution and thus picking a random move in OuterSearch might have a too small probability to yield a good move. Here our algorithm can be significantly improved by using better approaches to solve the mTSP or a more selective way to choose steps in OuterSearch. 

\begin{figure*}[ht] 
  \begin{subfigure}[b]{0.33\linewidth}
    \centering
    \includegraphics[width=0.75\linewidth]{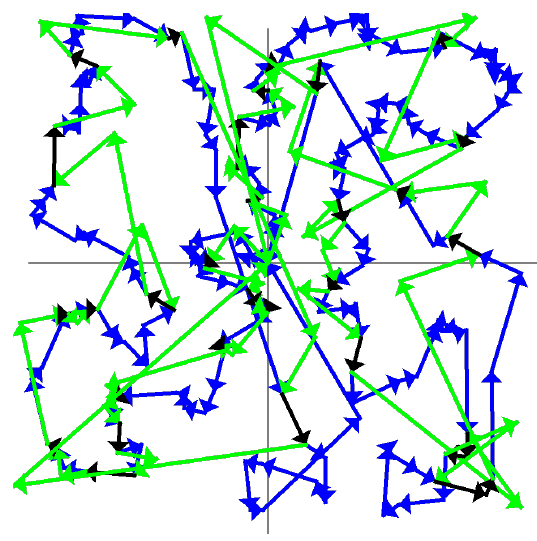} 
    \caption{All tours together}  
  \end{subfigure} 
  \begin{subfigure}[b]{0.33\linewidth}
    \centering
    \includegraphics[width=0.75\linewidth]{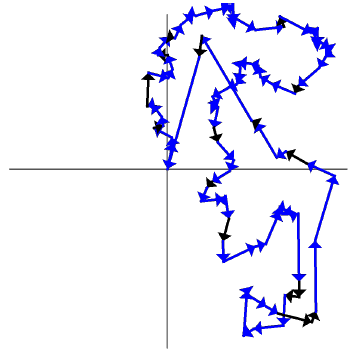} 
    \caption{Tour of first truck} 
  \end{subfigure} 
  \begin{subfigure}[b]{0.33\linewidth}
    \centering
    \includegraphics[width=0.75\linewidth]{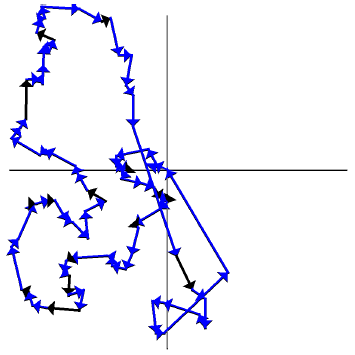} 
    \caption{Tour of second truck}  
  \end{subfigure}
  
  \begin{subfigure}[b]{0.33\linewidth}
    \centering
    \includegraphics[width=0.75\linewidth]{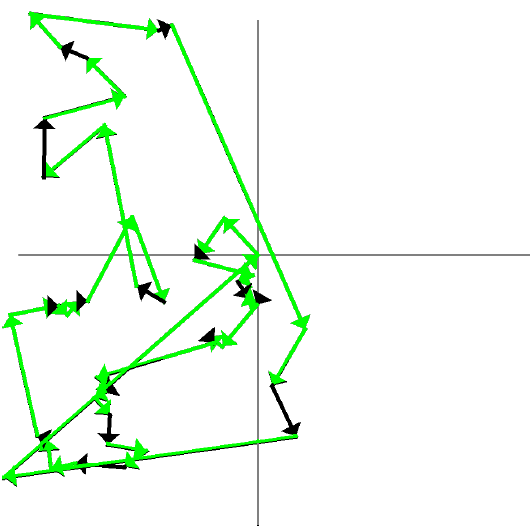} 
    \caption{Tour of first drone} 
  \end{subfigure} 
    \begin{subfigure}[b]{0.33\linewidth}
    \centering
    \includegraphics[width=0.75\linewidth]{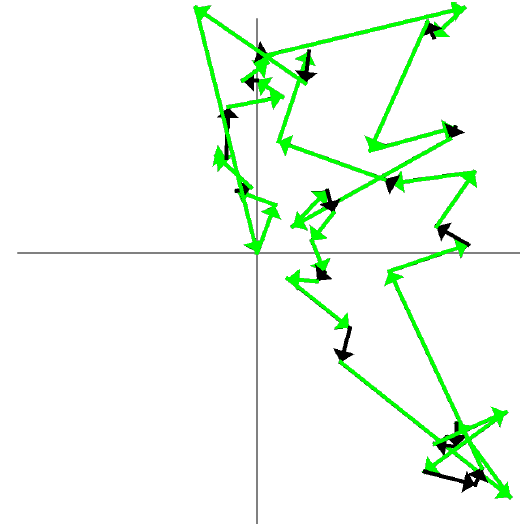} 
    \caption{Tour of second drone} 
  \end{subfigure}
  \caption{An example solution with two trucks and two drones}
  \label{exampleSolution} 
\end{figure*}

For an example solution see figure \ref{exampleSolution}. There we consider a sample from the computations for table \ref{parameters} number 2; in particular there are 200 packages, 2 trucks and 2 drones. 
Trucks carrying at least one drone are visualised in black, trucks not carrying any drone are visualised in blue, drones driving on a truck are also represented in black and flying drones are represented in green.
Its objective value is 972. Considering the solution two aspects are apparent: First, in the beginning, the distances between packages are much smaller than in the very end, which is good, as our objective function is the average delivery time of a package and thus long delivery times for a few packages in the end are not bad if in return we are able to deliver many packages very fast in the beginning.  Second we note, that the drones are used almost all the time.

We will now study the impacts of different parameters in more detail.

Therefore we will fix some of the parameters "number of packages", "number of trucks" and "number of drones" while varying the others. We do this according to table \ref{runningparameters}.

\begin{table}
\caption{}
\label{runningparameters}
 \begin{tcolorbox}[tab2,tabularx={p{0.5cm}p{3cm} p{3cm}},title=Settings for Running Parameters]
  \textbf{\#} & \textbf{running} & \textbf{fixed} \\ \hline
  1& \# packages (25 / 50 / 75 / 100 / 125 / 150 / 200 / 250 / 300) & 2 trucks, 2 drones  \\ \hline 
  2 & \# drones (1 / 2 / 3 / 4 / 5)& 2 trucks, 200 packages \\ \hline 
 3 & \# trucks = \# drones (1 / 2 / 3 / 4 / 5) & 200 packages 
\end{tcolorbox}
\end{table}

\begin{figure*}
\includegraphics[scale=0.7]{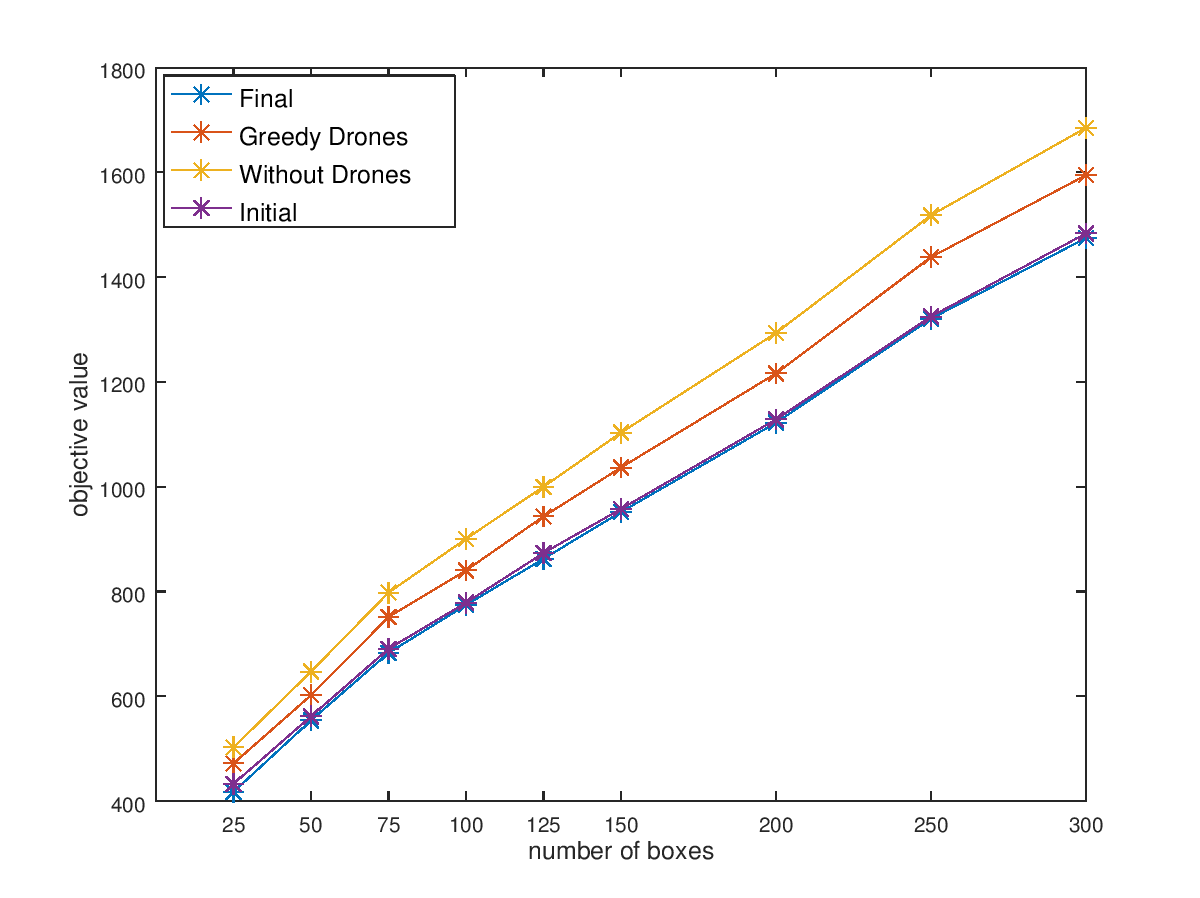}
\caption{Computational Results according to the data from table \ref{runningparameters} number 1.}
\label{tab4num1}
\end{figure*}

\begin{figure*}
\includegraphics[scale=0.7]{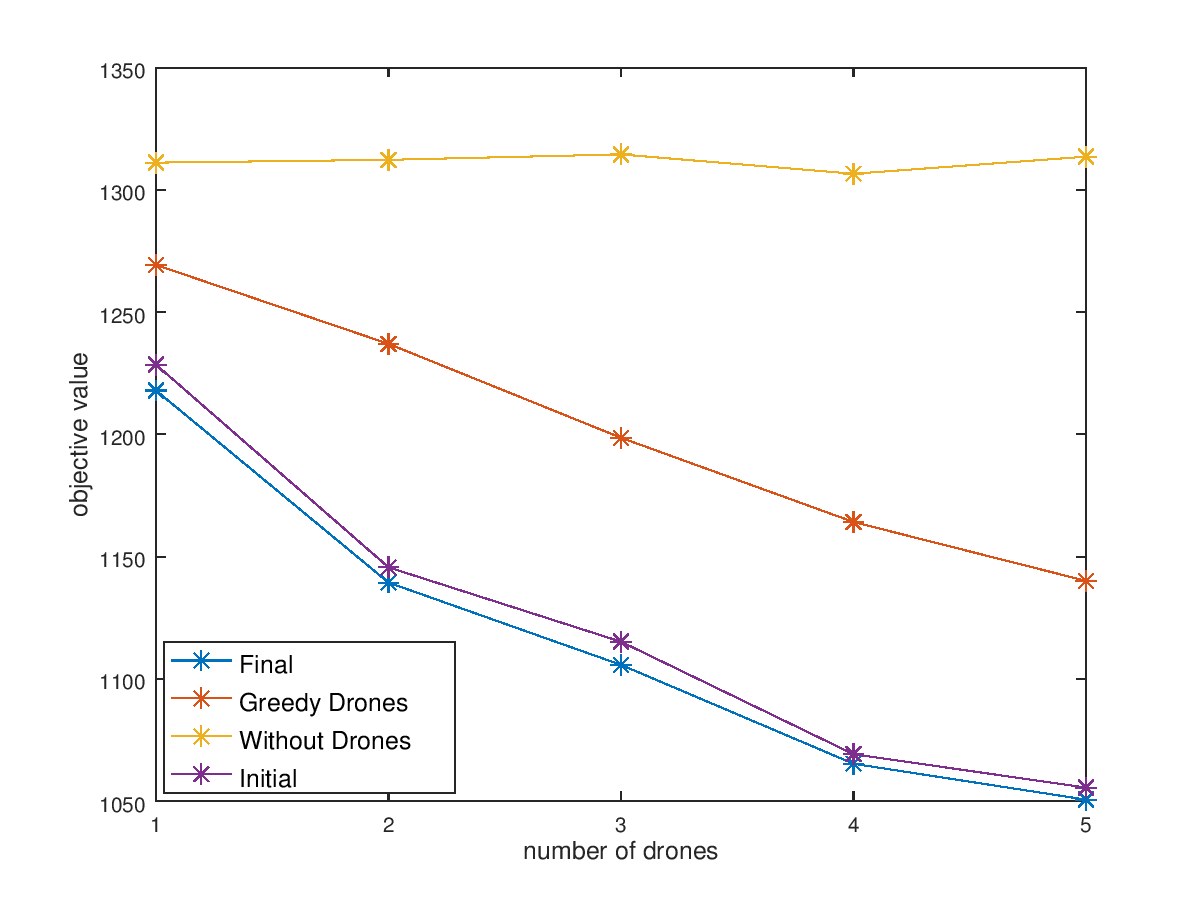}
\caption{Computational Results according to the data from table \ref{runningparameters} number 2.}
\label{tab4num2}
\end{figure*}

\begin{figure*}
\includegraphics[scale=0.7]{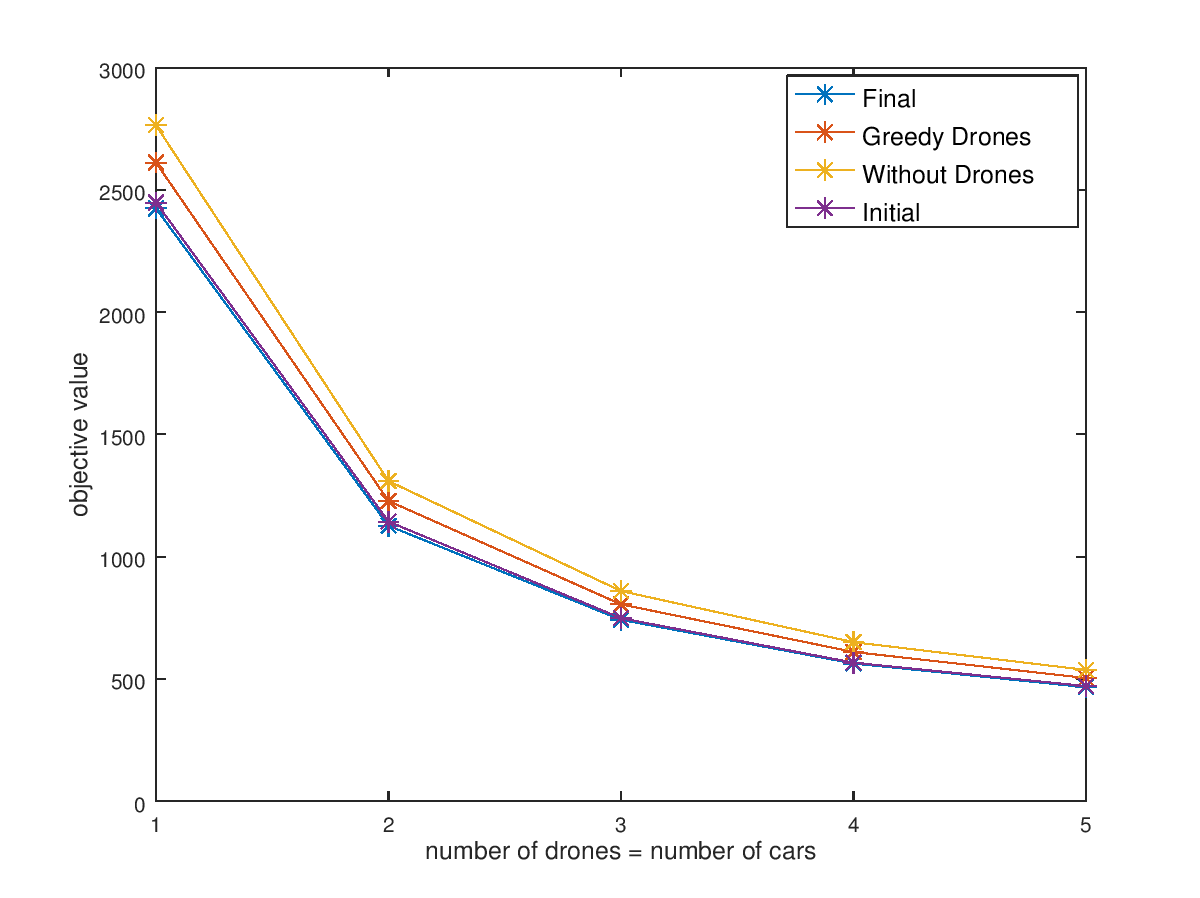}
\caption{Computational Results according to the data from table \ref{runningparameters} number 3.}
\label{tab4num3}
\end{figure*}

The results are visualised in figures \ref{tab4num1}, \ref{tab4num2} and \ref{tab4num3}; every setting from table 4 number 1 was averaged over 10 different instances, i.e. 10 different instances with 25 packages, 10 different instances with 50 packages etc. For table 4 number 2 every setting was sampled over the same 10 instances (for better comparability); the same holds for table 4 number 3.

So in all three scenarios Greedy Drones is outperformed significantly. If the number of packages is increased the objective value (recall: the objective value is the average delivery time of a package) appears to be slightly sublinear but close to linear.\footnote{Not in a strict sense but what the total appearance of the curve suggests.} The increase of the number of drones always improves the solution, however, the use of every new drone decreases.
The same qualitative behaviour occurs if the number of drones and the number of trucks is increased simultaneously.
\section{Conclusion and Outlook} \label{outlook}

We introduced the VRD problem and tackled it with nested local search algorithms using \cite{james}. Therefore we had to introduce suitable neighbourhood definitions. It turned out that our algorithm produces solutions which are significantly better than the Greedy Drones approach we introduced. However, the impact of OuterSearch was only marginal. An explanation has been discussed. Also we have seen that the use of drones improves the solution quality considerably suggesting that logistic processes may benefit from the use of drones. 

Note also that our approach can easily be adapted, e.g. further constraints can easily be added. Possible improvements for our algorithms can be achieved by further improving the neighbourhood definitions, choosing a more sophisticated probability distribution for picking neighbours as well as improving the method to obtain an solution for the mTSP.

This might also increase the impact on OuterSearch, as until now the individual truck tours might be too strongly separated to encourage changes like for example sending a drone from one truck to another. 

There are numerous further aspects that deserve attention:
\begin{itemize}
\item We did not consider the packing of the packages. It would be interesting to consider also packing constraints for the trucks. These could be comparatively simple ones like restricting the number of packages per truck or more complicated ones like $k$-dimensional bin packing, $k=1,2,3$.
\item A variant would be to allow the drones to land on the trucks even when they are driving. 
\item The starting and landing of a drone could cost some time. Also other simplifications we made could be dropped. 
\item A more realistic model of the charging process would be desirable, e.g. that drones charge with a certain speed, and can even leave when their battery is only charged partially. 
\item Basically all restrictions and variants of vehicle routing problems could also be applied here, e.g. time windows in which certain packages have to be delivered et cetera. 
\end{itemize}

\section{Acknowledgments} 
We would like to express our deep gratitude to our advisor Professor Martin Schottenloher for his valuable and constructive suggestions as well as his steady support.

\phantomsection
\bibliographystyle{abbrv}
\small
\bibliography{sample}

\setcounter{tocdepth}{0}
\normalsize

\addtocontents{toc}{\protect\setcounter{tocdepth}{2}}

\section{Appendix}
\subsection{Proof of Theorem \ref{almostfeasiblethm}}

\begin{proof}
We want to show an equivalence, so we have to prove the following two implications:
\begin{enumerate}
\item Every almost feasible solution which is inconsistent has a cycle without flip and not containing the depot in its graph.
\item 
If we have a cycle in the graph of a feasible solution, then it is a flip cycle or contains the depot.
\end{enumerate}
To show the first implication we assume that we have an almost feasible solution which violates the second condition of consistency. Thus there is a vehicle $V_1$ at some point in time which cannot continue, no matter how long it waits. So there is another vehicle $V_2$ it waits for, with which it is going to travel with. Again, $V_2$ waits for a vehicle, with which it is going to travel with and so on. As there are only finitely many vehicles, we obtain a chain $V_i \sim V_{i+1} \sim ... \sim V_n \sim V_i$ where "$A \sim B$" means "$A$ waits for $B$ to travel with" ($i$ may be $1$ but doesn't have to). We call the respective package positions where the vehicles got stuck $P_i,...,P_n$. Note that none of these positions can be the depot, as every vehicle that arrives at the depot stays there. 

As our solution is almost feasible, the respective vehicles would (if they would not have to wait) travel along edges from the solution, as follows:
$V_i$ travels to $P_n$ and travels on together with $V_n$. They then may stay together or split, but they travelled at least one edge together. Then $V_n$ travels to $P_{n-1}$ and then $V_n$ travels together at least one edge with $V_{n-1}$, and then $V_{n-1}$ continues to travel to $P_{n-2}$ and then $V_{n-1}$ travels at least one edge together with $V_{n-2}$ and so on. Finally $V_{i+1}$ travels to $P_{i}$ and travels at least one edge together with $V_i$ and then $V_i$ continues to travel to $P_n$ and they travel at least one edge together.

Now we can construct a cycle in the graph of the solution which contains no flips. Therefore we take edges along those edges on which vehicle $i$ would travel to vehicle $n$, those on which vehicle $n$ would travel to vehicle $n-1$ and so on. Again, none of the respective nodes can be the depot, as vehicles that reach the depot stay there. However, until now we only specified the start and end nodes of the edges that we choose, but there still are several possibilities to choose in our multi graph as between two nodes there may be several edges representing different vehicles. We choose always the edge which represents the vehicle travelling, e.g. for those edges on which vehicle $i$ travels to vehicle $n$ we choose the edges associated with vehicle $i$. However, there is one exception to that: Edges after any $P_j$ (for all $j$) 
are chosen as truck edges which is possible because two vehicles travel together there and thus one of them has to be a truck. This guarantees that we obtain a cycle which is not a flip cycle: Two consecutive edges represent either the same vehicle or one of them is an edge representing a truck. 
However until now we ignored one subtlety: We only constructed a circuit. To show that we can obtain a cycle we have to show that no node has to be visited more than once. We do this by either obtaining a contradiction or showing that a circuit that visits a node more than once can be replaced by one that visits it only once. Thus assume that there is a node that is visited more than once. This means that there are at least two ingoing and two outgoing edges. In particular the node must also be visited by exactly one truck and at least two of the four edges must represent a vehicle which either is a truck or travels with a truck. As there is only one truck per node allowed and this truck must perform a tour (in particular, no node is visited twice), it follows that exactly two of the four edges represent vehicles which either are trucks or represent vehicles that travel with trucks. It is clear that these two edges consist of one ingoing and one outgoing edge. Thus we are in the basic situation described in figure \ref{nocyclusbasic}.

 \begin{figure}
\includegraphics[scale=1.0]{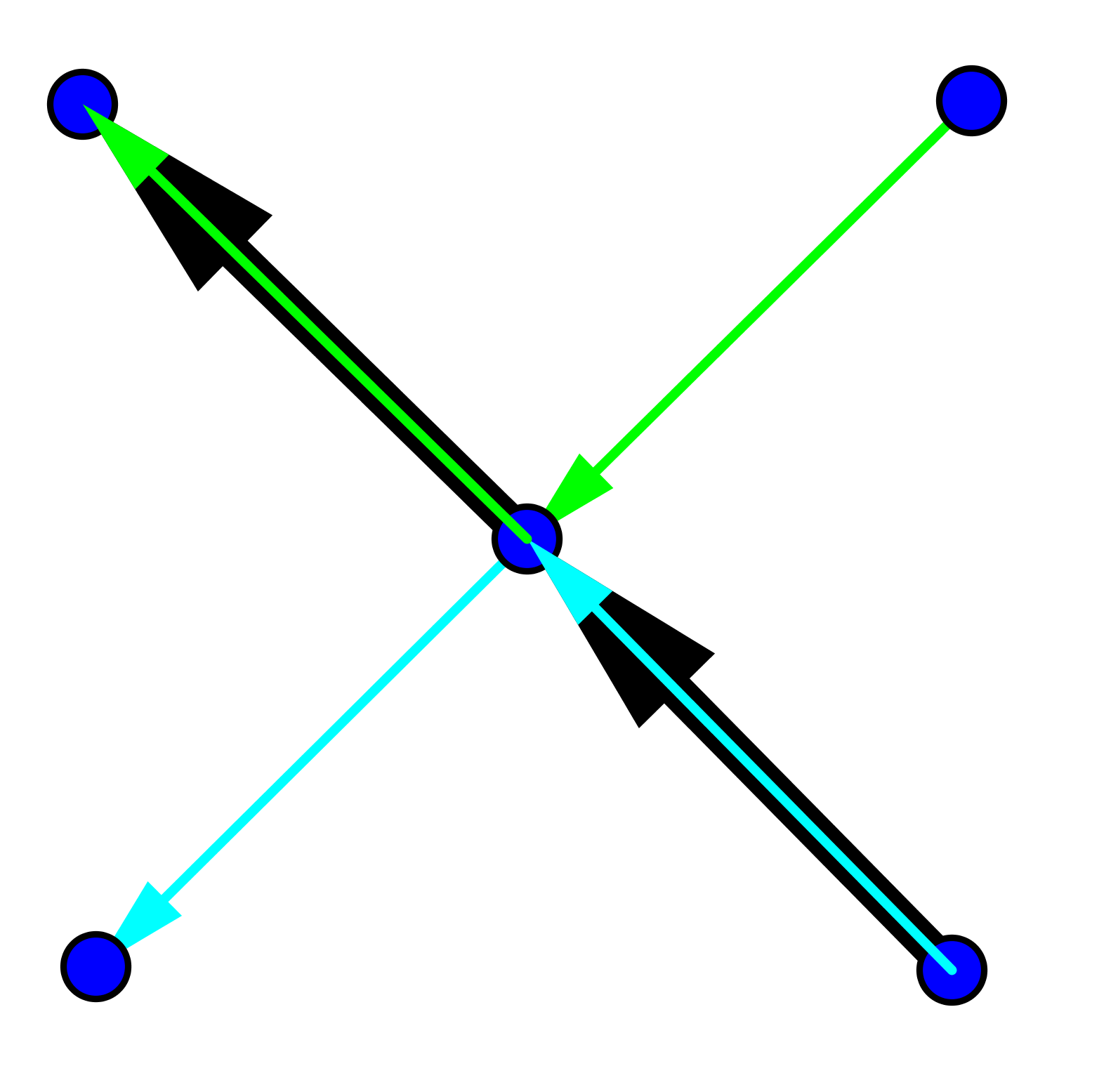}
\caption{Illustration to the argument about circles and circuits} 
\label{nocyclusbasic}
\end{figure}

However, as the constructed circuit has no flips, the green drone and the blue drone must be the same drone, but this contradicts the fact that the drone tours itself do not visit any node more than once in any almost feasible solution.\\

Now we show the converse direction. Therefore assume we have an almost feasible solution which has a cycle that is not a flip-cycle and does not contain the depot. We then have to show that it is not consistent.

For each node (by assumption) there are two possible scenarios how it is visited. Either by a drone only or by a truck, the drones that it carries when it arrives and the drones that it carries when it leaves.

Thus we have to distinguish different cases.
There are essentially 9 cases that we have to distinguish, which are illustrated in figure \ref{9cases}. Note that in figure \ref{9cases} green edges do not necessarily represent the same drone.
Also note again that the depot is not part of the cycle.
Clearly cases 6 and 8 cannot occur. Thus we look at the seven remaining cases. For each we show that edge $e$ cannot start before the previous edge in the cycle which gives us a contradiction, as the node under consideration cannot be the depot. 
Note that the cycle given to us is a priori given without the carry function or any information about which drone or truck is represented by a respective edge. Thus we actually only obtain four different possibilities for consecutive edges, either "drone-drone", "drone-truck", "truck-drone" or "truck-truck". However, each of these four possibilities corresponds to one of the nine cases illustrated, so if we excluded the nine cases we excluded in particular the four cases.

\begin{enumerate}
\item Cases 1,3,7,9: The incoming truck and the outgoing truck are the same truck, thus these cases are clear.
\item Case 2:
This case is clear, as only drones which arrived on the truck may leave the package position alone.
\item Case 4: Here the claim follows as the incoming drone must leave the package position again and can do this only by taking the truck (it cannot fly away, because then the presence of a truck would not have been allowed). As there can be at most one truck, it has to leave with that truck. So before $e$ is travelled, the previous edge has to be travelled.
\item Case 5: In this case we have to distinguish two subcases:
If the green edges belong to the same drone tour both times, then everything is clear as the respective drone is the only drone which visits the package position. However, the other subcase is that here a flip is represented. But we excluded this case by assumption (as otherwise the statement would have become wrong).
\end{enumerate}
 \begin{figure*}
\includegraphics[scale=1]{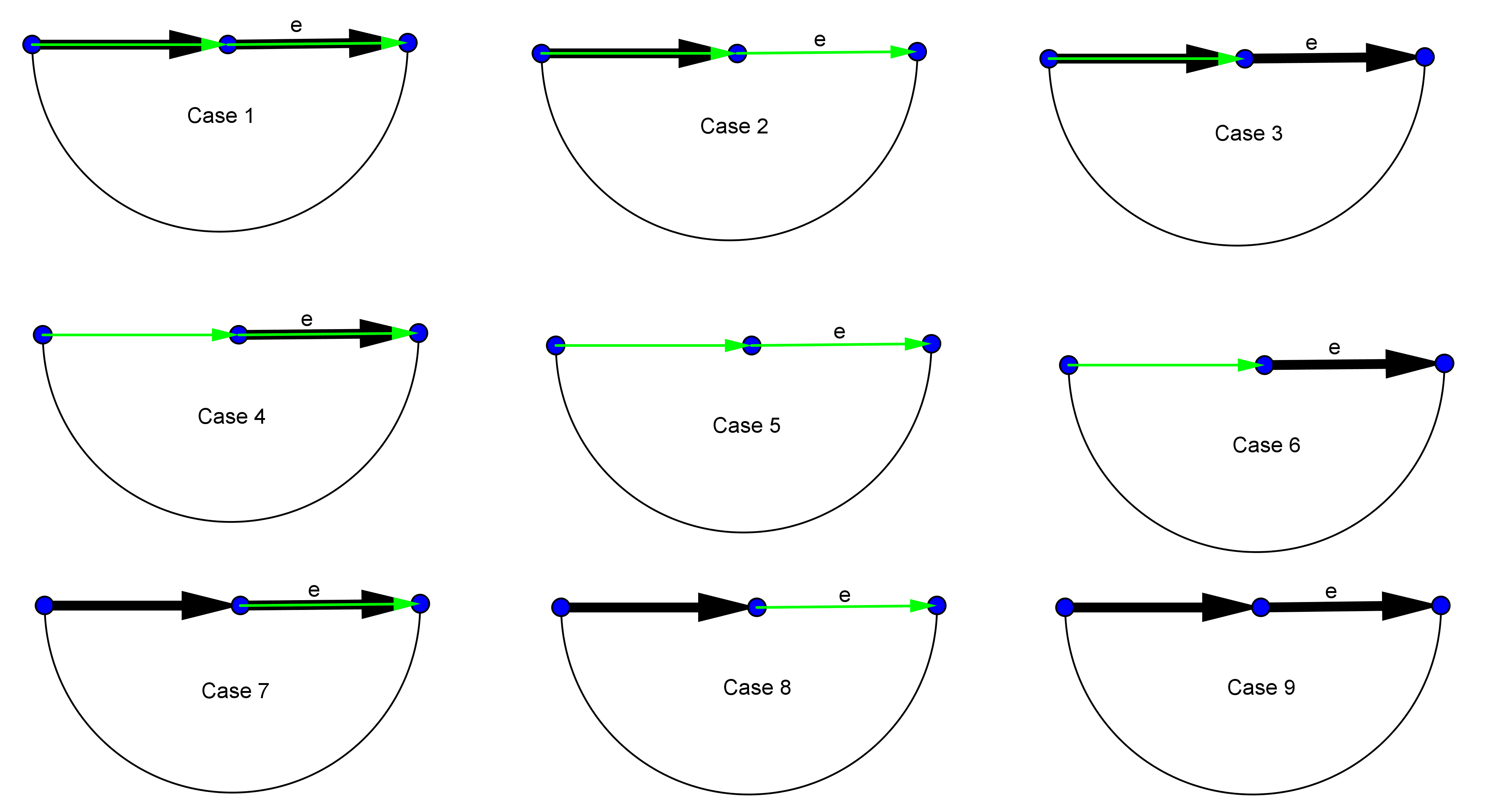}
\caption{9 cases we distinguish}
\label{9cases}
\end{figure*}
In particular we have a cycle where each edge cannot be performed before the previous one which is a contradiction to the feasibility. 
\end{proof}

\onecolumn

\subsection{Data of Table 3 --- Extensive Computations with Long Running Time} 

The settings for the TSP, SpeedUp1, OuterSearch and SpeedUp2 are the same for all 3 numbers of table 3.
All results are averaged about 10 samples; they are rounded to three decimals.

\subsubsection{Settings of TSP}

\begin{mdframed}
stopCriterion= [max steps without improvement: 1000000]

 jamesMethod=ParallelTempering
 
  minTemperature= 1.0E-6
  
   maxTemperature=10.0
   
numReplicasForParallelTempering=3
\end{mdframed}

\subsubsection{Settings of SpeedUp1 to be used for the initial solution}

\begin{mdframed}
stopCriterion= [max steps without improvement: 3000]

 jamesMethod=ParallelTempering
 
  minTemperature=0.001
  
   maxTemperature=1.0
   
numReplicasForParallelTempering=3
\end{mdframed}

\subsubsection{Settings of OuterSearch}

\begin{mdframed}
stopCriterion= [max steps without improvement: 5]

 jamesMethod=ParallelTempering
 
  minTemperature=1.0E-7
  
   maxTemperature=10.0
   
numReplicasForParallelTempering=3
\end{mdframed}

\subsubsection{Settings of SpeedUp2 to be used in the OuterSearch}

\begin{mdframed}[nobreak=true]
stopCriterion= [max steps without improvement: 5]

 jamesMethod=ParallelTempering
 
  minTemperature=1.0E-6

   maxTemperature=10.0
   
numReplicasForParallelTempering=3
\end{mdframed}

\subsubsection{Table 3 Number 1}

\paragraph{Basic Settings}

\begin{mdframed}

numberOfStartAnglesForTsp: 10

packageRange: 200

numberOfPackages:200

numberOfTrucks:1

numberOfDrones:2

reachOfDrones: 10000

numberOfSamples:10
\end{mdframed}

\paragraph{Results for Table 3 Number 1}

\begin{center}
    \begin{tabular}{|p{6cm} | p{6cm} |}
    \hline
    Final                  & 2095.624 \\ \hline   
   Initial                 &  2105.923 \\ \hline
   Without Drones          &  2542.843 \\ \hline
   Greedy Drones                 &  2307.287 \\ \hline
    \end{tabular}
\end{center}

Running time (in sec.) for the each of the 10 samples and in average over these samples:

\begin{center}
    \begin{tabular}{|p{1.2cm} | p{1.2cm} |p{1.2cm} |p{1.2cm} |p{1.2cm} |p{1.2cm} |p{1.2cm} |p{1.2cm} |p{1.2cm} |p{1.2cm} |p{1.2cm} |}
    \hline
   \# 1 & \# 2& \# 3&\# 4&\# 5&\# 6&\# 7&\# 8&\# 9&\# 10& $\varnothing$ \\ \hline
    32694 & 28372  & 28252 & 27321 & 25188 & 32230 & 24777 & 28052 & 25671 & 27618 & 28017.5\\ \hline
    \end{tabular}
\end{center}

\subsubsection{Table 3 Number 2}

\paragraph{Basic Settings}

\begin{mdframed}

numberOfStartAnglesForTsp: 10

packageRange: 200

numberOfPackages: 200

numberOfTrucks: 2

numberOfDrones: 2

reachOfDrones: 10000

numberOfSamples: 10
\end{mdframed}

\paragraph{Results for Table 3 Number 2}

\begin{center}
    \begin{tabular}{|p{6cm} | p{6cm} |}
    \hline
    Final                   & 1070.495 \\ \hline   
   Initial                  & 1082.144   \\ \hline
   Without Drones           & 1243.052  \\ \hline
   Greedy Drones                  & 1164.001  \\ \hline
    \end{tabular}
\end{center}

Running time (in sec.) for the each of the 10 samples and in average over these samples:

\begin{center}
    \begin{tabular}{|p{1.2cm} | p{1.2cm} |p{1.2cm} |p{1.2cm} |p{1.2cm} |p{1.2cm} |p{1.2cm} |p{1.2cm} |p{1.2cm} |p{1.2cm} |p{1.2cm} |}
    \hline
   \# 1 & \# 2& \# 3&\# 4&\# 5&\# 6&\# 7&\# 8&\# 9&\# 10& $\varnothing$ \\ \hline
     37712 & 50628  & 43936 & 37158 & 38688 & 40243 & 39695 & 37909 & 38295 & 39142 & 40340.6\\ \hline
    \end{tabular}
\end{center}

\subsubsection{Table 3 Number 3}

\paragraph{Basic Settings}

\begin{mdframed}

numberOfStartAnglesForTsp: 10

packageRange: 200

numberOfPackages: 200

numberOfTrucks: 3

numberOfDrones: 5

reachOfDrones: 10000

numberOfSamples:10
\end{mdframed}

\paragraph{Results for Table 3 Number 3}

\begin{center}
    \begin{tabular}{|p{6cm} | p{6cm} |}
    \hline
    Final                   & 688.603 \\ \hline   
   Initial                  & 700.780  \\ \hline
   Without Drones           & 833.028 \\ \hline
   Greedy Drones                  & 758.635 \\ \hline
    \end{tabular}
\end{center}

Running time (in sec.) for the each of the 10 samples and in average over these samples:

\begin{center}
    \begin{tabular}{|p{1.2cm} | p{1.2cm} |p{1.2cm} |p{1.2cm} |p{1.2cm} |p{1.2cm} |p{1.2cm} |p{1.2cm} |p{1.2cm} |p{1.2cm} |p{1.2cm} |}
    \hline
   \# 1 & \# 2& \# 3&\# 4&\# 5&\# 6&\# 7&\# 8&\# 9&\# 10& $\varnothing$ \\ \hline
    62301 & 61520  & 72879 & 69580 & 59586 & 71220 & 68349 & 61023 & 77505 & 70641 & 67460.4 \\ \hline
    \end{tabular}
\end{center}

\subsection{Data of Table 4}
The settings for the TSP, SpeedUp1, OuterSearch and SpeedUp2 are the same for all 3 numbers of table 4.
All results are averaged about 10 samples; they are rounded to three decimals.

\subsubsection{Basic Settings for TSP}

\begin{mdframed}
stopCriterion=[max runtime: 55000 ms]

 jamesMethod=ParallelTempering
 
  minTemperature=1.0E-6
  
   maxTemperature=10.0
   
    numReplicasForParallelTempering=3
\end{mdframed}

\subsubsection{Settings of SpeedUp1 to be used for the initial solution}

\begin{mdframed}
stopCriterion=[max steps without improvement: 1750]

 jamesMethod=ParallelTempering
 
  minTemperature=0.001
  
   maxTemperature=1.0
   
numReplicasForParallelTempering=3
\end{mdframed}

\subsubsection{Basic Settings for OuterSearch}

\begin{mdframed}
stopCriterion=[max steps without improvement: 2]

 jamesMethod=ParallelTempering
 
  minTemperature=1.0E-7
  
   maxTemperature=10.0

     numReplicasForParallelTempering=3
     \end{mdframed}

\subsubsection{Basic Settings for SpeedUp2 to be used in OuterSearch}

\begin{mdframed}
stopCriterion=[max steps without improvement: 2]

 jamesMethod=ParallelTempering
 
  minTemperature=1.0E-6
  
   maxTemperature=10.0

     numReplicasForParallelTempering=3
     \end{mdframed}
     
\subsubsection{Data of Table 4 Number 1 --- Number of Packages is Varying}

\paragraph{Basic Settings}

\begin{mdframed}
numberOfStartAnglesForTsp: 5

packageRange: 200

numberOfTrucks: 2

numberOfDrones:2

reachOfDrones: 10000

numberOfSamples: 10
\end{mdframed}

\paragraph{Results for numberOfPackages = 25}

\begin{center}
    \begin{tabular}{|p{6cm} | p{6cm} |}
    \hline
    Final                   & 417,587 \\ \hline   
   Initial                  & 432,700 \\ \hline
   Without Drones           & 502,596 \\ \hline
   Greedy Drones                  & 471,827 \\ \hline
    \end{tabular}
\end{center}

Running time (in sec.) for the each of the 10 samples and in average over these samples:

\begin{center}
    \begin{tabular}{|p{1.2cm} | p{1.2cm} |p{1.2cm} |p{1.2cm} |p{1.2cm} |p{1.2cm} |p{1.2cm} |p{1.2cm} |p{1.2cm} |p{1.2cm} |p{1.2cm} |}
    \hline
   \# 1 & \# 2& \# 3&\# 4&\# 5&\# 6&\# 7&\# 8&\# 9&\# 10& $\varnothing$ \\ \hline
   1655 & 1636 & 1619 & 1573 & 1584 & 1648 &1597 &1638 &1633 &1594 & 1617.7\\ \hline
    \end{tabular}
\end{center}

\paragraph{Results for numberOfPackages = 50}

\begin{center}
    \begin{tabular}{|p{6cm} | p{6cm} |}
    \hline
 Final                   & 554.137 \\ \hline   
   Initial                  & 561.772 \\ \hline
   Without Drones           & 647.094 \\ \hline
   Greedy Drones                  & 602.057 \\ \hline
    \end{tabular}
\end{center}

Running time (in sec.) for the each of the 10 samples and in average over these samples:

\begin{center}
    \begin{tabular}{|p{1.2cm} | p{1.2cm} |p{1.2cm} |p{1.2cm} |p{1.2cm} |p{1.2cm} |p{1.2cm} |p{1.2cm} |p{1.2cm} |p{1.2cm} |p{1.2cm} |}
    \hline
   \# 1 & \# 2& \# 3&\# 4&\# 5&\# 6&\# 7&\# 8&\# 9&\# 10& $\varnothing$ \\ \hline
    2314 &2135  & 2326 & 2238 &2291  & 2083 & 2371& 2144&2362 & 2068& 2233.2 \\ \hline
    \end{tabular}
\end{center}

\paragraph{Results for numberOfPackages = 75}

\begin{center}
    \begin{tabular}{|p{6cm} | p{6cm} |}
    \hline
   Final                   & 682.728 \\ \hline   
   Initial                  & 690.603 \\ \hline
   Without Drones           & 797.851 \\ \hline
   Greedy Drones                  & 751.164 \\ \hline
    \end{tabular}
\end{center}

Running time (in sec.) for the each of the 10 samples and in average over these samples:

\begin{center}
    \begin{tabular}{|p{1.2cm} | p{1.2cm} |p{1.2cm} |p{1.2cm} |p{1.2cm} |p{1.2cm} |p{1.2cm} |p{1.2cm} |p{1.2cm} |p{1.2cm} |p{1.2cm} |}
    \hline
   \# 1 & \# 2& \# 3&\# 4&\# 5&\# 6&\# 7&\# 8&\# 9&\# 10& $\varnothing$ \\ \hline
    2486& 2931 & 2676 & 2746 & 2693  & 2936 &2675 &3174 & 2590& 2925 & 2783.2 \\ \hline
    \end{tabular}
\end{center}

\paragraph{Results for numberOfPackages = 100}

\begin{center}
    \begin{tabular}{|p{6cm} | p{6cm} |}
    \hline
    Final                   & 774.205 \\ \hline   
   Initial                  & 778.659 \\ \hline
   Without Drones           & 900.642 \\ \hline
   Greedy Drones                  & 840.174 \\ \hline
    \end{tabular}
\end{center}

Running time (in sec.) for the each of the 10 samples and in average over these samples:

\begin{center}
    \begin{tabular}{|p{1.2cm} | p{1.2cm} |p{1.2cm} |p{1.2cm} |p{1.2cm} |p{1.2cm} |p{1.2cm} |p{1.2cm} |p{1.2cm} |p{1.2cm} |p{1.2cm} |}
    \hline
   \# 1 & \# 2& \# 3&\# 4&\# 5&\# 6&\# 7&\# 8&\# 9&\# 10& $\varnothing$ \\ \hline
    3567& 3460  & 3946 & 3739 & 3279 &3540  &3174 &3530 & 3135&3204 & 3457.4\\ \hline
    \end{tabular}
\end{center}

\paragraph{Results for numberOfPackages = 125}

\begin{center}
    \begin{tabular}{|p{6cm} | p{6cm} |}
    \hline
   Final                   & 862.451  \\ \hline   
   Initial                  & 873.390  \\ \hline
   Without Drones           & 999.880  \\ \hline
   Greedy Drones                  & 943.257  \\ \hline
    \end{tabular}
\end{center}

Running time (in sec.) for the each of the 10 samples and in average over these samples:

\begin{center}
    \begin{tabular}{|p{1.2cm} | p{1.2cm} |p{1.2cm} |p{1.2cm} |p{1.2cm} |p{1.2cm} |p{1.2cm} |p{1.2cm} |p{1.2cm} |p{1.2cm} |p{1.2cm} |}
    \hline
   \# 1 & \# 2& \# 3&\# 4&\# 5&\# 6&\# 7&\# 8&\# 9&\# 10& $\varnothing$ \\ \hline
    5496 & 4004 & 4291 & 4475 & 3986 & 5322 & 3953 & 3525&4873 & 4166 & 4409.1\\ \hline
    \end{tabular}
\end{center}

\paragraph{Results for numberOfPackages = 150}

\begin{center}
    \begin{tabular}{|p{6cm} | p{6cm} |}
    \hline
    Final                    & 951.777  \\ \hline   
   Initial                   & 957.611  \\ \hline
   Without Drones            & 1103.267 \\ \hline
   Greedy Drones                   & 1036.633 \\ \hline
    \end{tabular}
\end{center}

Running time (in sec.) for the each of the 10 samples and in average over these samples:

\begin{center}
    \begin{tabular}{|p{1.2cm} | p{1.2cm} |p{1.2cm} |p{1.2cm} |p{1.2cm} |p{1.2cm} |p{1.2cm} |p{1.2cm} |p{1.2cm} |p{1.2cm} |p{1.2cm} |}
    \hline
   \# 1 & \# 2& \# 3&\# 4&\# 5&\# 6&\# 7&\# 8&\# 9&\# 10& $\varnothing$ \\ \hline
     4030 & 4299 & 4909 & 4886 & 4638 & 4854 &5288 &3930 &4287 & 5025 &4614.6 \\ \hline
    \end{tabular}
\end{center}

\paragraph{Results for numberOfPackages = 200}

\begin{center}
    \begin{tabular}{|p{6cm} | p{6cm} |}
    \hline
   Final                    & 1122.241  \\ \hline   
   Initial                   & 1128.471  \\ \hline
   Without Drones            & 1293.464 \\ \hline
   Greedy Drones                   & 1216.198 \\ \hline
    \end{tabular}
\end{center}

Running time (in sec.) for the each of the 10 samples and in average over these samples:

\begin{center}
    \begin{tabular}{|p{1.2cm} | p{1.2cm} |p{1.2cm} |p{1.2cm} |p{1.2cm} |p{1.2cm} |p{1.2cm} |p{1.2cm} |p{1.2cm} |p{1.2cm} |p{1.2cm} |}
    \hline
   \# 1 & \# 2& \# 3&\# 4&\# 5&\# 6&\# 7&\# 8&\# 9&\# 10& $\varnothing$ \\ \hline
   6509 & 5802 & 6926 & 5978 & 6022 & 8137 & 7730& 6922 &7045 & 5030 &6610.1 \\ \hline
    \end{tabular}
\end{center}

\paragraph{Results for numberOfPackages = 250}

\begin{center}
    \begin{tabular}{|p{6cm} | p{6cm} |}
    \hline
    Final                    & 1320.723 \\ \hline   
   Initial                   & 1324.807 \\ \hline
   Without Drones            & 1518.106 \\ \hline
   Greedy Drones                   & 1438.265 \\ \hline
    \end{tabular}
\end{center}

Running time (in sec.) for the each of the 10 samples and in average over these samples:

\begin{center}
    \begin{tabular}{|p{1.2cm} | p{1.2cm} |p{1.2cm} |p{1.2cm} |p{1.2cm} |p{1.2cm} |p{1.2cm} |p{1.2cm} |p{1.2cm} |p{1.2cm} |p{1.2cm} |}
    \hline
   \# 1 & \# 2& \# 3&\# 4&\# 5&\# 6&\# 7&\# 8&\# 9&\# 10& $\varnothing$ \\ \hline
    6085 & 7108 & 6886 & 7261 & 6460 & 8152 & 7890 & 6636 & 12055 & 6454 & 7498.7\\ \hline
    \end{tabular}
\end{center}

\paragraph{Results for numberOfPackages = 300}

\begin{center}
    \begin{tabular}{|p{6cm} | p{6cm} |}
    \hline
   Final                   & 1474.680 \\ \hline   
   Initial                  & 1483.482 \\ \hline
   Without Drones           & 1685.132 \\ \hline
   Greedy Drones                  & 1594.704 \\ \hline
    \end{tabular}
\end{center}

Running time (in sec.) for the each of the 10 samples and in average over these samples:

\begin{center}
    \begin{tabular}{|p{1.2cm} | p{1.2cm} |p{1.2cm} |p{1.2cm} |p{1.2cm} |p{1.2cm} |p{1.2cm} |p{1.2cm} |p{1.2cm} |p{1.2cm} |p{1.2cm} |}
    \hline
   \# 1 & \# 2& \# 3&\# 4&\# 5&\# 6&\# 7&\# 8&\# 9&\# 10& $\varnothing$ \\ \hline
   7407 & 9034 & 8144 & 8215 & 8751 & 8113 & 7874 & 18973 & 8342& 10116 & 9496.9\\ \hline
    \end{tabular}
\end{center}

\subsubsection{Data of Table 4 Number 2 --- Number of Drones is Varying}

\paragraph{Basic Settings}

\begin{mdframed}

numberOfStartAnglesForTsp: 5

packageRange: 200

numberOfTrucks: 2

numberOfPackages: 200

reachOfDrones: 10000

numberOfSamples: 10
\end{mdframed}

Note that the 10 samples are the same for all numbers of drones.

\paragraph{Results for 1 Drone}

\begin{center}
    \begin{tabular}{|p{6cm} | p{6cm} |}
    \hline
    Final                    & 1217.935 \\ \hline   
   Initial                   & 1228.510 \\ \hline
   Without Drones            & 1311.129 \\ \hline
   Greedy Drones                    & 1269.300 \\ \hline
    \end{tabular}
\end{center}

Running time (in sec.) for the each of the 10 samples and in average over these samples:

\begin{center}
    \begin{tabular}{|p{1.2cm} | p{1.2cm} |p{1.2cm} |p{1.2cm} |p{1.2cm} |p{1.2cm} |p{1.2cm} |p{1.2cm} |p{1.2cm} |p{1.2cm} |p{1.2cm} |}
    \hline
   \# 1 & \# 2& \# 3&\# 4&\# 5&\# 6&\# 7&\# 8&\# 9&\# 10& $\varnothing$ \\ \hline
    3581 & 6856 & 3894 & 3728 & 4746 & 6425 & 3582 & 3936 & 3412 & 3897 & 4405.7 \\ \hline
    \end{tabular}
\end{center}

\paragraph{Results for 2 Drones}

\begin{center}
    \begin{tabular}{|p{6cm} | p{6cm} |}
    \hline
   Final                   & 1139.378 \\ \hline   
   Initial                  & 1145.609 \\ \hline
   Without Drones           & 1312.271\\ \hline
   Greedy Drones                  & 1237.011 \\ \hline
    \end{tabular}
\end{center}

Running time (in sec.) for the each of the 10 samples and in average over these samples:

\begin{center}
    \begin{tabular}{|p{1.2cm} | p{1.2cm} |p{1.2cm} |p{1.2cm} |p{1.2cm} |p{1.2cm} |p{1.2cm} |p{1.2cm} |p{1.2cm} |p{1.2cm} |p{1.2cm} |}
    \hline
   \# 1 & \# 2& \# 3&\# 4&\# 5&\# 6&\# 7&\# 8&\# 9&\# 10& $\varnothing$ \\ \hline
    6181 & 5079 & 9478 & 5423 & 5275 & 6513 & 5875 & 7935 & 5229 & 5363 & 6235.1 \\ \hline
    \end{tabular}
\end{center}

\paragraph{Results for 3 Drones}

\begin{center}
    \begin{tabular}{|p{6cm} | p{6cm} |}
    \hline
    Final                    & 1105.702 \\ \hline   
   Initial                   & 1115.183 \\ \hline
   Without Drones            & 1314.491 \\ \hline
   Greedy Drones                   & 1198.601 \\ \hline
    \end{tabular}
\end{center}

Running time (in sec.) for the each of the 10 samples and in average over these samples:

\begin{center}
    \begin{tabular}{|p{1.2cm} | p{1.2cm} |p{1.2cm} |p{1.2cm} |p{1.2cm} |p{1.2cm} |p{1.2cm} |p{1.2cm} |p{1.2cm} |p{1.2cm} |p{1.2cm} |}
    \hline
   \# 1 & \# 2& \# 3&\# 4&\# 5&\# 6&\# 7&\# 8&\# 9&\# 10& $\varnothing$ \\ \hline
    8934 & 8255 & 8219 & 12658 & 8725 & 7611 & 10911 & 8523 & 8216 & 6749 & 8880.1\\ \hline
    \end{tabular}
\end{center}

\paragraph{Results for 4 Drones}

\begin{center}
    \begin{tabular}{|p{6cm} | p{6cm} |}
    \hline
    Final                    & 1065.267 \\ \hline   
   Initial                   & 1069.172\\ \hline
   Without Drones            & 1306.574\\ \hline
   Greedy Drones                   & 1164.155 \\ \hline
    \end{tabular}
\end{center}

Running time (in sec.) for the each of the 10 samples and in average over these samples:

\begin{center}
    \begin{tabular}{|p{1.2cm} | p{1.2cm} |p{1.2cm} |p{1.2cm} |p{1.2cm} |p{1.2cm} |p{1.2cm} |p{1.2cm} |p{1.2cm} |p{1.2cm} |p{1.2cm} |}
    \hline
   \# 1 & \# 2& \# 3&\# 4&\# 5&\# 6&\# 7&\# 8&\# 9&\# 10& $\varnothing$ \\ \hline
   12385 & 11151 & 10982 & 9152 & 10432 & 12924 & 13484 & 10673 & 9430 & 10774 & 11138.7\\ \hline
    \end{tabular}
\end{center}

\paragraph{Results for 5 Drones}

\begin{center}
    \begin{tabular}{|p{6cm} | p{6cm} |}
    \hline
    Final                    & 1050.583  \\ \hline   
   Initial                   & 1055.523  \\ \hline
   Without Drones            & 1313.569 \\ \hline
   Greedy Drones                   & 1140.065 \\ \hline
    \end{tabular}
\end{center}

Running time (in sec.) for the each of the 10 samples and in average over these samples:

\begin{center}
    \begin{tabular}{|p{1.2cm} | p{1.2cm} |p{1.2cm} |p{1.2cm} |p{1.2cm} |p{1.2cm} |p{1.2cm} |p{1.2cm} |p{1.2cm} |p{1.2cm} |p{1.2cm} |}
    \hline
   \# 1 & \# 2& \# 3&\# 4&\# 5&\# 6&\# 7&\# 8&\# 9&\# 10& $\varnothing$ \\ \hline
   12985 & 12875  & 12373 & 14029 & 14512 & 13057 & 12482 & 14987 & 13293 & 12240 & 13283.3\\ \hline
    \end{tabular}
\end{center}

\subsubsection{Data of Table 4 Number 3 --- Number of Drones = Number of Trucks is Varying}

\paragraph{Basic Settings}

\begin{mdframed}

numberOfStartAnglesForTsp: 5

packageRange: 200

numberOfPackages: 200

reachOfDrones: 10000

numberOfSamples: 10
\end{mdframed}

Note that the 10 samples are the same for all numbers of drones (= numbers of trucks).

\paragraph{Results for number of trucks = number of drones = 1}

\begin{center}
    \begin{tabular}{|p{6cm} | p{6cm} |}
    \hline
    Final                    & 2423.913  \\ \hline   
   Initial                   & 2448.041  \\ \hline
   Without Drones            & 2764.246  \\ \hline
   Greedy Drones                   & 2611.550  \\ \hline
    \end{tabular}
\end{center}

Running time (in sec.) for the each of the 10 samples and in average over these samples:

\begin{center}
    \begin{tabular}{|p{1.2cm} | p{1.2cm} |p{1.2cm} |p{1.2cm} |p{1.2cm} |p{1.2cm} |p{1.2cm} |p{1.2cm} |p{1.2cm} |p{1.2cm} |p{1.2cm} |}
    \hline
   \# 1 & \# 2& \# 3&\# 4&\# 5&\# 6&\# 7&\# 8&\# 9&\# 10& $\varnothing$ \\ \hline
   5304 & 5769  & 4759 & 4761 & 5137 & 5393 & 6337 & 4514 & 4524 & 4098 & 5059.6 \\ \hline
    \end{tabular}
\end{center}

\paragraph{Results for number of trucks = number of drones = 2}

\begin{center}
    \begin{tabular}{|p{6cm} | p{6cm} |}
    \hline
    Final                   & 1126.035 \\ \hline   
   Initial                  & 1143.066  \\ \hline
   Without Drones           & 1308.897 \\ \hline
   Greedy Drones                  & 1227.583 \\ \hline
    \end{tabular}
\end{center}

Running time (in sec.) for the each of the 10 samples and in average over these samples:

\begin{center}
    \begin{tabular}{|p{1.2cm} | p{1.2cm} |p{1.2cm} |p{1.2cm} |p{1.2cm} |p{1.2cm} |p{1.2cm} |p{1.2cm} |p{1.2cm} |p{1.2cm} |p{1.2cm} |}
    \hline
   \# 1 & \# 2& \# 3&\# 4&\# 5&\# 6&\# 7&\# 8&\# 9&\# 10& $\varnothing$ \\ \hline
    6340 & 6652  & 5859 & 7496 & 7167 & 8419 & 7068 & 6599 & 7201 & 6229 & 6903 \\ \hline
    \end{tabular}
\end{center}

\paragraph{Results for number of trucks = number of drones = 3}

\begin{center}
    \begin{tabular}{|p{6cm} | p{6cm} |}
    \hline
   Final                 & 741.483 \\ \hline   
   Initial                & 748.102  \\ \hline
   Without Drones         & 859.068 \\ \hline
   Greedy Drones                & 804.655 \\ \hline
    \end{tabular}
\end{center}

Running time (in sec.) for the each of the 10 samples and in average over these samples:

\begin{center}
    \begin{tabular}{|p{1.2cm} | p{1.2cm} |p{1.2cm} |p{1.2cm} |p{1.2cm} |p{1.2cm} |p{1.2cm} |p{1.2cm} |p{1.2cm} |p{1.2cm} |p{1.2cm} |}
    \hline
   \# 1 & \# 2& \# 3&\# 4&\# 5&\# 6&\# 7&\# 8&\# 9&\# 10& $\varnothing$ \\ \hline
    12955 & 8030  & 10463 & 8375 & 8141 & 8565 & 8019 & 9112 & 10165 & 6917 & 9074.2\\ \hline
    \end{tabular}
\end{center}

\paragraph{Results for number of trucks = number of drones = 4}

\begin{center}
    \begin{tabular}{|p{6cm} | p{6cm} |}
    \hline
     Final                   & 562.574 \\ \hline   
   Initial                  & 566.079  \\ \hline
   Without Drones           & 649.992 \\ \hline
   Greedy Drones                  & 609.802 \\ \hline
    \end{tabular}
\end{center}

Running time (in sec.) for the each of the 10 samples and in average over these samples:

\begin{center}
    \begin{tabular}{|p{1.2cm} | p{1.2cm} |p{1.2cm} |p{1.2cm} |p{1.2cm} |p{1.2cm} |p{1.2cm} |p{1.2cm} |p{1.2cm} |p{1.2cm} |p{1.2cm} |}
    \hline
   \# 1 & \# 2& \# 3&\# 4&\# 5&\# 6&\# 7&\# 8&\# 9&\# 10& $\varnothing$ \\ \hline
   9241 & 10131  & 10910 & 10731 & 10868 & 9205 & 14559 & 15193 & 10400 & 14379 & 11561.7\\ \hline
    \end{tabular}
\end{center}

\paragraph{Results for number of trucks = number of drones = 5}

\begin{center}
    \begin{tabular}{|p{6cm} | p{6cm} |}
    \hline
    Final                   & 465.627  \\ \hline   
   Initial                  & 470.237  \\ \hline
   Without Drones           & 536.365 \\ \hline
   Greedy Drones                   & 503.363 \\ \hline
    \end{tabular}
\end{center}

Running time (in sec.) for the each of the 10 samples and in average over these samples:

\begin{center}
    \begin{tabular}{|p{1.2cm} | p{1.2cm} |p{1.2cm} |p{1.2cm} |p{1.2cm} |p{1.2cm} |p{1.2cm} |p{1.2cm} |p{1.2cm} |p{1.2cm} |p{1.2cm} |}
    \hline
   \# 1 & \# 2& \# 3&\# 4&\# 5&\# 6&\# 7&\# 8&\# 9&\# 10& $\varnothing$ \\ \hline
    16217 & 13302  & 11438 & 14861 & 13040 & 11587 & 13150 & 12427 & 12557 & 12610 & 13118.9\\ \hline
    \end{tabular}
\end{center}


\end{document}